\documentclass[12pt]{article}
\usepackage{latexsym}
\usepackage{amsmath,amsfonts,amssymb,epsfig}
\usepackage{amsthm,amscd}
\usepackage{mathrsfs} 
\usepackage{hyperref}
\usepackage{cancel}
\usepackage{appendix}

\numberwithin{equation}{section}

\newtheorem{theorem}{Theorem}[section]
\newtheorem{definition}[theorem]{Definition}

\newtheorem{proposition}[theorem]{Proposition}
\newtheorem{lemma}[theorem]{Lemma}

\newtheorem{rem}[theorem]{Remark}
\newenvironment{remark}{\begin{rem} \rm}{\end{rem}}
\newtheorem{ide}[theorem]{Idea}

\newtheorem{exa}[theorem]{Example}

\newtheorem{ques}[theorem]{Question}

\newtheorem{pro}[theorem]{Problem}

\newtheorem{spe}[theorem]{Speculation}

\font\bbb=msbm10 scaled 1100

\newcommand{\Lie}{{\cal L}}
\newcommand{\R}{\mbox{\bbb R}}       

\newcommand{\Z}{\mbox{\bbb Z}}
\newcommand{\Cnf}{\text{\rm Conf}}

\newcommand{\no}{\noindent}

\title{The third order helicity of magnetic fields\\ via link maps.} 

\author{R. Komendarczyk\footnote{\emph{2000 Mathematics Subject Classification}. Primary: 76W05,57M25, Secondary: 58F18,58A10
}\ \footnote{
  This project is supported  by DARPA, \#FA9550-08-1-0386.}
}

\begin{document}
\maketitle

\begin{abstract}
 \no We introduce an alternative approach to the third order helicity of a volume preserving vector field $B$, which leads us to a lower bound for the $L^2$-energy of $B$. The proposed approach exploits correspondence between the Milnor $\bar{\mu}_{123}$-invariant for 3-component links and the homotopy invariants of maps to configuration spaces, and we provide a simple geometric proof of this fact in the case of Borromean links. 
 Based on these connections we develop a formulation for the third order helicity of $B$ on invariant \emph{unlinked} domains of $B$, and 
 provide Arnold's style ergodic interpretation of this invariant
 as an average asymptotic $\bar{\mu}_{123}$-invariant of orbits of $B$. 
\end{abstract}


\section{Introduction}
 A purpose of this paper is to develop a particular formula for the third order helicity on certain  invariant sets of a volume preserving vector field $B$. The third order helicity, \cite{Khesin98}, is an invariant of $B$ under the action of volumorphisms isotopic to the identity (denoted here by $\text{\rm SDiff}_0(M)$).  
\no Importance of such invariants stems from the basic fact that the evolution of the vorticity in the ideal hydrodynamics or of the magnetic field $B_0$ in the ideal magnetohydrodynamics (MHD), occurs along a path $t\longrightarrow g(t)\in \text{\rm SDiff}_0(M)$, \cite[p. 176]{Khesin98}. Namely,
$B(t)=g_\ast(t)B_0$ which is a direct consequence of Euler's equations: 
\begin{equation}\label{eq:Eulers-eq}
\frac{d}{dt} B+[v, B]=0,\qquad \frac{d}{dt} g(t)=v .
\end{equation}
 One often says that the magnetic field $B$ is \emph{frozen in} the velocity field $v$ of plasma, and the action by $\text{\rm SDiff}_0(M)$ is frequently referred to as \emph{frozen-in-field} deformations. A fundamental example of such invariant defined for a general class in $\text{\rm SVect}(M)$ is the \emph{helicity} $\mathsf{H}_{12}(B_1,B_2)$ defined for a pair of vector fields $B_1$ and $B_2$ on $M=S^3$ or a homology sphere. 
Helicity has been first introduced by Woltjer, \cite{Woltjer58}, in the context of magnetic fields, and is a measure of how orbits of $B_1$, and $B_2$ link with each other. This topological interpretation of helicity has been made precise by Arnold, who introduced the concept of the average asymptotic linking number of  a volume preserving vector field $B$ on $M$, \cite{Arnold86}. 
The subject has been further investigated
in \cite{Akhmetiev05, Berger90, Laurence-Stredulinsky00b, Hornig04} (see \cite{Khesin98} for additional references), where authors approach \emph{higher helicities} via the Massey products under various assumptions about the vector fields or their domains (work in \cite{Gambaudo-Ghys97, Spera06, Verjovsky94} concerns yet other approaches to the problem). Extensions of the helicity concept to higher dimensional 
foliations can be found in \cite{Khesin92, Riviere02, Kotschick-Vogel03} and recently in
\cite{Cantarella-Parsley09}.

In this paper we present an alternative to these approaches, which is a natural extension of the notion of the linking number as a degree of a map, and exploits relations to the homotopy theory of certain maps associated to the link. The paper has two parts.

In the first part we show how the Milnor $\bar{\mu}_{123}$-invariant: $\bar{\mu}_{123}(L)$ of a parametrized Borromean link $L$ in $S^3$ can be obtained as a Hopf degree of an associated map to the configuration space of three points in $S^3$. The presented approach to link homotopy invariants of 3-component links has been proposed in \cite{Kohno02} for $n$-component links in $\R^3$, and conveniently simplified for 3-component links in $S^3$ in the joint work \cite{Deturck-Gluck-Melvin-Shonkwiler08}, where the full correspondence between $\bar{\mu}_{123}(L)$ and the Pontryagin-Hopf degree is proved. In Section \ref{sec:mu-and-hopf} we present the original proof of this correspondence in the Borromean case, which is sufficient for our purposes.  

The second part of the paper is discusses a new definition of the third order helicity denoted here by $\mathsf{H}_{123}(B;\mathcal{T})$, where $B$ is a volume preserving vector field having an invariant unlinked domain $\mathcal{T}\subset S^3$. The simplest of such domains are three invariant 
handlebodies in $S^3$ which have pairwise unlinked 
cycles in the first homology. This includes the case of Borromean flux tubes already investigated in  \cite{Berger90, Laurence-Stredulinsky00b}. In Section \ref{sec:ergodic}  we develop an ergodic formulation of $\mathsf{H}_{123}(B;\mathcal{T})$ as an average asymptotic $\bar{\mu}_{123}$-invariant, in the spirit of Arnold's average asymptotic linking number, which allows us to extend the definition of the invariant to topologically more complicated unlinked domains. We also derive, in Section \ref{sec:energy},  a lower 
bound for the $L^2$-energy of $B$ in terms of $\mathsf{H}_{123}(B;\mathcal{T})$.

\vspace{1cm}

\no\emph{Acknowledgements:} The inspiration for the presented approach to the link homotopy invariants  comes from the paper of Toshitake Kohno \cite{Kohno02}, and I am grateful to him for the valuable e-mail correspondence. I have enjoyed conversations with many colleagues at the University of Pennsylvania, who have influenced this work. I wish to thank Herman Gluck for weekly meetings and his interest in this project, Frederic Cohen, Dennis DeTurck, Charlie Epstein, Paul Melvin, Tristan Rivi$\grave{\rm e}$re, Clay Shonkwiler, Jim Stasheff, David Shea Vick for the valuable input.  I am also grateful to my advisor 
Robert Ghrist who introduced me to the subject long time ago.  After posting recent joint results in \cite{Deturck-Gluck-Melvin-Shonkwiler08} we were informed by Paul Kirk about related works of Urlich Koschorke in \cite{Koschorke97, Koschorke04}, on homotopy invariants of link maps and Milnor $\bar{\mu}$-invariants. The author acknowledges financial support of DARPA, \#FA9550-08-1-0386.


\section{The Milnor \texorpdfstring{$\bar{\mu}_{123}$}{mu-123}-invariant and the Hopf degree.}\label{sec:mu-and-hopf}

The $\bar{\mu}$-invariants of $n$-component links in $S^3$ have been introduced by Milnor in \cite{Milnor54, Milnor57}
as invariants of links up to \emph{link homotopy}. Recall that the link homotopy 
is a deformation of a link in $S^3$ which allows each component to pass through 
itself but not through a different component. Clearly, this is a weaker equivalence than the equivalence of links up to \emph{isotopy} where components are not allowed to pass through themselves at all. The fundamental example of a $\bar{\mu}$-invariant is the linking number (denoted by $\bar{\mu}_{12}$) which is a complete invariant of the 2-component links up to link homotopy. In the realm of 3-component links the relevant invariants are the pairwise linking numbers $\bar{\mu}_{12}$, $\bar{\mu}_{23}$, $\bar{\mu}_{32}$, and the third invariant $\bar{\mu}_{123}$ in $\mathbb{Z}_{\text{gcd}(\bar{\mu}_{12}, \bar{\mu}_{23}, \bar{\mu}_{32})}$, which is 
a well defined integer, if and only if, $\bar{\mu}_{12}=\bar{\mu}_{23}=\bar{\mu}_{32}=0$. In the second part of the paper we will interpret this statement as a topological condition on the invariant set of a vector  field. A precise definition of $\bar{\mu}$-invariants is algebraic and involves the Magnus expansion of the lower central series of the fundamental group: $\pi_1(S^3-L)$ of the link complement. We refer the interested reader to the works in \cite{Milnor54, Milnor57}. In the remaining part of this section we will prove that $\bar{\mu}_{123}(L)$ is a Hopf degree for an appropriate map associated to the link $L$, provided that 
the link is \emph{Borromean}, i.e. the pairwise linking numbers are zero (note that the Borromean links are more general then  Brunnian links, \cite{Milnor54}).

Let us review basic facts about the Hopf degree $\mathscr{H}(f)$ of a map 
$f:S^3\longrightarrow S^2$, (see
e.g. \cite{Bott82}). A well known property of the Hopf degree is that  $\mathscr{H}: f\longrightarrow \mathscr{H}(f)$ provides an isomorphism between $\pi_3(S^2)$
and $\Z$. Recall that up to a constant multiple we
may express $\mathscr{H}(f)$ as ($M=S^3$)
\begin{gather}\label{eq:hopf-integral}
\mathscr{H}(f)= \int_{M} \alpha\wedge f^\ast\nu=\int_{M} \alpha\wedge \omega=\int_{M} \alpha\wedge d\alpha,
\end{gather}
where $\nu$ is the area 2-form on $S^2$, and $\alpha$ satisfies
$\omega=f^\ast\nu=d\alpha$. Notice that $f^\ast\nu$ is always exact
since the cohomology of $S^3$ in dimension $2$ vanishes. We may also interpret $\mathscr{H}(f)$ as an intersection number,
\cite{Bott82, Milnor97}. Namely, consider two regular values $p_1$ and $p_2\in S^2$ of the map $f$,
then $l_1=f^{-1}(p_1)$ and $l_2=f^{-1}(p_2)$ form a link in $S^3$, and the integral formula
\eqref{eq:hopf-integral} can be interpreted, as the intersection number of $l_1$ with the Seifert surface spanning $l_2$:
\begin{gather}\label{eq:hopf-intersection}
\mathscr{H}(f)=\text{lk}(l_1,l_2)\ .
\end{gather}
If we replace $S^3$ with an arbitrary closed
compact orientable 3-dimensional manifold $M$, we may still obtain
an invariant of $f:M\longrightarrow S^2$ this way, provided that
the condition $f^\ast\nu=d\alpha$ holds. 

\begin{proposition}\label{th:hopf-degree-general}
 Let $M$ be a closed Riemannian manifold, and $\nu\in \Omega^2(S^2)$ the area form on $S^2$.
 The formula \eqref{eq:hopf-integral} provides a homotopy invariant 
 for a map $f:M\longrightarrow S^2$, if the 2-form $f^\ast\nu$ is exact. Up to a constant multiple 
 this invariant can be calculated as an intersection number defined in \eqref{eq:hopf-intersection},
 where $l_1=f^{-1}(p_1)$ and $l_2=f^{-1}(p_2)$ form a link in $M$, where both $l_1$ and $l_2$ are null-homologous.
\end{proposition}

\begin{proof}
\no Given a homotopy $F:I\times M\mapsto
S^2$, $f_1=F(1,\,\cdot\,)$, $f_0=F(0,\,\cdot\,)$, we define
$\hat{\omega}=F^\ast\nu$. We have 
\[
\hat{\omega}  = F^\ast\nu=d\,\hat{\alpha},\qquad 
\omega_1  =  f_1^\ast\nu=i^\ast_1\,F^\ast\nu=d\,\alpha_1,\qquad 
\omega_0  =  f_0^\ast\nu=i^\ast_0\,F^\ast\nu=d\,\alpha_0,
\]
where $i_0:M\hookrightarrow M\times I$, $i_0(x)=(x,0)$, $i_1:M\hookrightarrow M\times I$, $i_1(x)=(x,1)$
are appropriate inclusions. Potentials: $\hat{\alpha}\bigl|_{\{0\}\times M}$, $\alpha_0$, and $\hat{\alpha}\bigl|_{\{1\}\times M}$, $\alpha_1$ differ  by a closed form 
\[
\hat{\alpha}\bigl|_{\{0\}\times M}-\alpha_0=\beta_0,\quad \hat{\alpha}\bigl|_{\{1\}\times M}-\alpha_1=\beta_1,\quad d\beta_0=d\beta_1=0,
\]
therefore the Stokes Theorem immediately implies that Formula \eqref{eq:hopf-integral} is independent of the choice 
of the potential. For the proof of invariance under homotopies we revoke the standard argument in \cite[p. 228]{Bott82}
\[
\begin{split}
 0 & = \int_{F(M\times I)} \nu\wedge\nu
 =\int_{M\times I} \hat{\omega}\wedge \hat{\omega} =  \int_{M\times I} \hat{\omega}\wedge d\,\hat{\alpha}\\
 & = \int_{M\times I} d(\hat{\omega}\wedge \hat{\alpha}) = \int_{\partial(M\times I)} \hat{\omega}\wedge
\hat{\alpha}  =  \int_{M} \omega_1\wedge \alpha_1- \int_{M}
\omega_0\wedge \alpha_0\\
& =  \int_{M} d\alpha_1 \wedge \alpha_1-\int_{M}
d\alpha_0\wedge \alpha_0 =  \mathscr{H}(f_1)-\mathscr{H}(f_0).
\end{split}
\]
\no The interpretation of $\mathscr{H}(f)$ as the intersection number \eqref{eq:hopf-intersection} is the same as in \cite[p. 230]{Bott82}.
\end{proof}

\no Given a 3-component parametrized link $L=\{L_1,L_2,L_3\}$ in $S^3$
we wish to associate a certain map $F_L:S^1\times S^1\times S^1\mapsto S^2$ to it, and interpret its Hopf degree as the Milnor $\bar{\mu}_{123}$-invariant. Recall the definition of the configuration space of $k$ points in $M$:
\begin{gather*}
 \Cnf_k(M):=\{(x_1,x_2,\ldots, x_k)\in (M)^k\,|\,x_i\neq x_j,\text{for}\ i\neq j\}.
\end{gather*}
As an introduction to the method we review the Gauss formula for the linking number of 
a 2-component link $L=\{L_1, L_2\}$ in $\R^3$. Denote parameterizations 
of components by  $L_1=\{x(s)\}$, $L_2=\{y(t)\}$ and consider the map
\begin{gather*}
 F_L:S^1\times S^1\stackrel{L}{\longrightarrow} \Cnf_2(\R^3) \stackrel{r}{\longrightarrow} S^2,\qquad L(s,t)=(x(s), y(t))\ .
\end{gather*}
where $r(x,y)=\frac{x-y}{\|x-y\|}$ is a retraction of $\Cnf_2(\R^3)$ onto $S^2$.
It yields the classical Gauss linking number formula:
\begin{gather*}
 \bar{\mu}_{12}(L)=\text{lk}(L_1,L_2)=\text{deg}(F_L),\qquad\quad \text{deg}(F_L)=\int_{S^1\times S^1} F^\ast_L(\nu),
\end{gather*}
where $\nu\in \Omega^2(S^2)$ is the 
area form on $S^2$. Consequently, the linking number $\text{lk}(L_1,L_2)$,
also known as the Milnor $\bar{\mu}_{12}$-invariant, can be obtained as the homotopy invariant of the map $F_L$ associated to $L$. 
Observe that homotopy classes $[S^1\times S^1,S^2]$ are 
isomorphic to $\Z$ and $\text{deg}:F\rightarrow \text{deg}(F)$ provides the isomorphism,
we also point out that as sets: $[S^k,\Cnf_n(\R^3)]=\pi_k(\Cnf_n(R^3))$, and $[S^k,\Cnf_n(S^3)]=\pi_k(\Cnf_n(S^3))$ (c.f. \cite[p. 421]{Hatcher02}). Thus considering the based homotopies, in context of the link homotopy of Borromean links, and base point free homotopies is equivalent in this setting.

In \cite{Kohno02, Koschorke97},
authors consider a natural extension of this approach to $n$-component parame-trized links $L$ in $\R^3$ by considering maps $F_L:S^1\times\ldots\times S^1\longrightarrow \Cnf_n(\R^3)$ and their homotopy classes, we refer to this type of maps loosely as \emph{link maps}, (c.f. \cite{Koschorke97}). In particular, Kohno \cite{Kohno02} proposed specific representatives of cohomology classes of the based loop space of $\Cnf_n(\R^3)$ as candidates for appropriate link homotopy invariants of $L$. It has been observed, in \cite{Deturck-Gluck-Melvin-Shonkwiler08}, that in the $3$-component case it is  beneficial to consider $\Cnf_3(S^3)$, and $L\subset S^3$, since the topology of $\Cnf_3(S^3)$ simplifies dramatically (in comparison to $\Cnf_3(\R^3)$). We review this simplification in the following paragraph as it is essential for the proof of the main theorem in this section.


\begin{figure}
\centering
\begin{picture}(380, 150)
  \put(0,0){\includegraphics[width=5.2in]{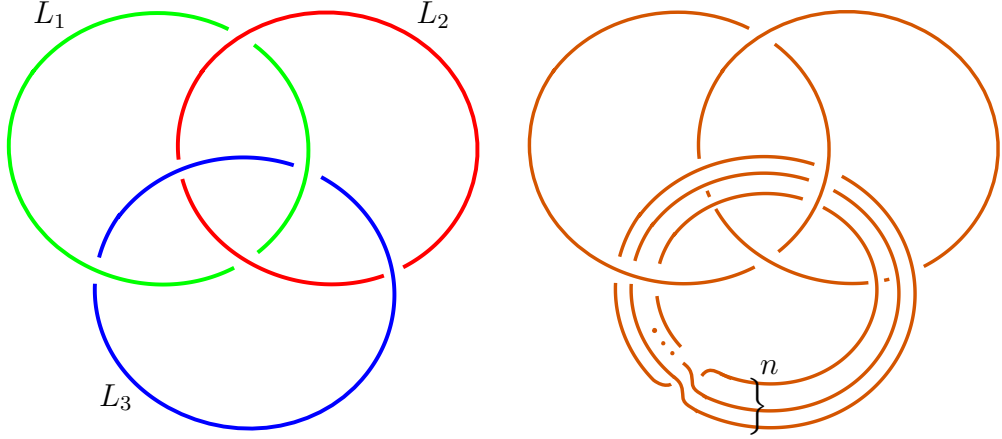}}
  \put(10,155){$L_1$} 
  \put(155,155){$L_2$} 
  \put(35,10){$L_3$} 
  \put(280,7){$\mathbf{\Bigl\}}$}
  \put(285,22){$n$}
\end{picture}  
\caption{left: $\bar{\mu}_{123}=\pm 1$ and right: $\bar{\mu}_{123}=\pm n$.}
\label{fig:n-borromean}
\end{figure}

Consider a 3-component link $L=\{L_1, L_2, L_3\}$ in $S^3$ 
parametrized by $\{x(s),y(t),z(u)\}$ 
and the following map 
\begin{gather}\label{eq:F}
 F_L:S^1\times S^1\times S^1\stackrel{L}{\longrightarrow} \Cnf_3(S^3)\stackrel{H}{\longrightarrow} S^2,\qquad L(s,t,u)=(x(s), y(t), z(u)),
\end{gather}
where we denote by 
$H:S^3\times \R^3\times (\R^3\setminus\{0\})\mapsto S^2$ the projection on the $S^2$ factor, 
first concluding that  
$\Cnf_3(S^3)\subset S^3\times S^3\times S^3$ is
diffeomorphic to $S^3\times \Cnf_2(\R^3)=S^3\times
\R^3\times (\R^3\setminus\{0\})$, and consequently deformation
retracts onto $S^3\times S^2$. 
Considering $S^3$ as unit quaternions, the map $H$ can be expressed explicitly by the formula, \cite{Deturck-Gluck-Melvin-Shonkwiler08}:
\begin{gather}\label{eq:hermans-map}
 \Cnf_3(S^3) \ni (x,y,z) \stackrel{H}{\longrightarrow}
 \frac{\text{pr}(x^{-1}\cdot y)-\text{pr}(x^{-1}\cdot z)}{\|\text{pr}(x^{-1}\cdot y)-\text{pr}(x^{-1}\cdot z)\|}\in S^2,
\end{gather}
where $\cdot$ stands for the quaternionic multiplication, $\ ^{-1}$
is the quaternionic inverse, and $\text{pr}:S^3\longrightarrow \R^3$ the
stereographic projection from $1$. As a result one has the following particular expression for $F_L$:
\begin{gather}\label{eq:F_L}
F_L(s,t,u)=\frac{\text{pr}(x(s)^{-1}\cdot z(u))-\text{pr}(x(s)^{-1}\cdot y(t))}{\|\text{pr}(x(s)^{-1}\cdot z(u))-\text{pr}(x(s)^{-1}\cdot
 y(t))\|}.
\end{gather}

 At this point we note that one has a freedom in choosing the deformation retraction $H$ in \eqref{eq:hermans-map}, but the above particular formula makes the proof of the main theorem of this section possible. Let $\mathbb{T}=S^1\times S^1\times S^1$ denote the domain of $F_L$, notice that, thanks to \eqref{eq:F_L}, restricting $F_L$ to the subtorus $\mathbb{T}_{23}$ in the second and third coordinate $(t,u)$ of $\mathbb{T}$, we obtain the usual Gauss map of the 2-component link $\{x^{-1}\cdot L_2,x^{-1} \cdot L_3\}$. Since the diffeomorphism $x^{-1}\cdot$ of $S^3$ is orientation preserving we conclude
 $\text{\rm deg}(F_L|_{\mathbb{T}_{23}})=\text{lk}(L_2,L_3)$. We claim that for any 2-component sublink 
 $\{L_i,L_j\}$ of $L$:
 \begin{equation}\label{eq:linking-numbers}
  \text{\rm deg}(F_L|_{\mathbb{T}_{ij}})=\pm\text{lk}(L_i,L_j),\qquad 1\leq i<j\leq 3,
 \end{equation}
 where $i,j$ index the coordinates of $\mathbb{T}$. Indeed, since already true for $i=2$ and $j=3$,  the general case follows by applying a permutation $\sigma\in \Sigma_3$ of coordinate factors in $\Cnf_3(S^3)\subset (S^3)^3$. Notice that $\sigma$ is a diffeomorphism of $\Cnf_3(S^3)\subset (S^3)^3$ either preserving or reversing the orientation (which explains the sign in \eqref{eq:linking-numbers}). We infer \eqref{eq:linking-numbers} because $\sigma$ induces an isomorphism on homotopy groups of $\Cnf_3(S^3)$.


\begin{figure}[htpb]
\begin{center}
\begin{picture}(414,140)
  \put(30,0){\includegraphics{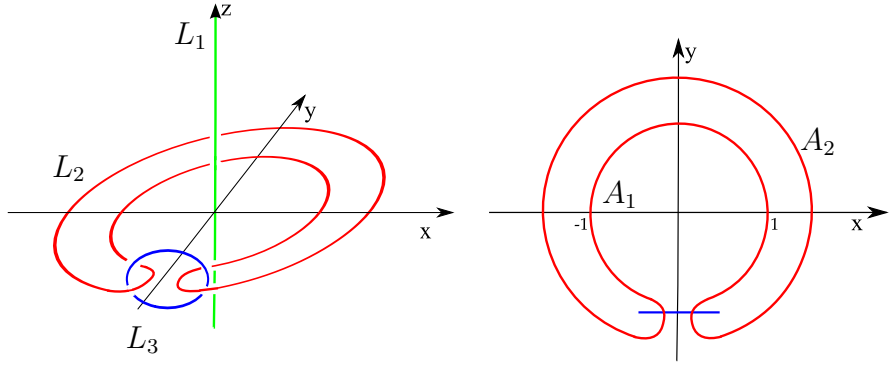}}
  \put(93,120){$L_1$}
  \put(47,70){$L_2$}
  \put(75,5){$L_3$}
  \put(330,80){$A_2$}
  \put(255,60){$A_1$}
\end{picture}  
\caption{The model of $L_\text{Borr}$ parametrized by $\{L_1=\text{pr}(x(s)),L_2=\text{pr}(y(t)),L_3=\text{pr}(z(u))\}$. Arc $A_1$ is a part of the unit circle on $xy$-plane, $A_2$ is a part of circle of radius $1+\epsilon$.} \label{fig:conf1}

\end{center}
\end{figure}

\no The main theorem of this section is
\begin{theorem}\label{th:milnor-hopf}
 Let $L=\{L_1,L_2,L_3\}$ be a 3-component Borromean link in $S^3$, consider
 the associated map $F_L$ defined in \eqref{eq:F}. The Hopf degree of this 
 map satisfies
 \begin{gather}
   \mathscr{H}(F_L)=\pm 2\,\bar{\mu}_{123},
 \end{gather}
 where the sign depends on the choice of orientations of components of $L$.
\end{theorem}

\begin{proof}

 By the \emph{Borromean rings} we understand any $3$-component link with $\bar{\mu}_{123}=\pm 1$. Every such link is link homotopic to the diagram presented on Figure \ref{fig:n-borromean} (where the sign can be determined from the orientation of components). The proof of the theorem can be reduced to the case of Borromean rings $L^\text{Borr}$, as follows:  if $L$ and $L'$ are link-homotopic then $F_{L}$ and $F_{L'}$ are homotopic maps.
 By the Milnor classification of 3-component links up to link homotopy (see \cite{Milnor54}), every 3-component link $L$ with zero pairwise linking numbers and $\bar{\mu}_{123}=\pm n$ is represented by the right diagram on Figure \ref{fig:n-borromean}. Consequently, up to homotopy, the associated map $F_L$ can be obtained from $F_{L_\text{Borr}}$ by covering one of the $S^1$ factors in $\mathbb{T}$, $n$-times. 
Therefore, in order to prove the claim it suffices to show 
\begin{equation}\label{eq:hopf=2}
 \mathscr{H}(F_{L^\text{Borr}})=\pm 2\ .
\end{equation}
 According to Proposition \ref{th:hopf-degree-general}, $\mathscr{H}(F_L)$ is well defined for a link $L\subset S^3$ provided  $F_L^\ast\nu\in \Omega^2(\mathbb{T})$ is trivial in $H^2(\mathbb{T})$, which is true thanks to \eqref{eq:linking-numbers} and because the pairwise linking numbers of $L$ are zero. The method of proof relies on a direct calculation of 
$\mathscr{H}(F_{L^\text{Bor}})$, for a carefully chosen parametrization of $L^\text{Bor}$ in $S^3$. This calculation is achieved by 
visualization of the link $l_{S,N}=l_S\cup l_N$ in $\mathbb{T}$, and application of Formula \eqref{eq:hopf-intersection}, where
\[
 l_S:=F^{-1}_{L^\text{Bor}}(S),\qquad l_N:=F^{-1}_{L^\text{Bor}}(N)
\]
are preimages of the North pole $N=(0,0,1)$ and South pole $S=(0,0,-1)$ in $S^2\subset \R^3$. Notice that $[F_L^\ast\nu]=0$ in $H^2(\mathbb{T})$, if and only if, $[l_S]=0$ and $[l_N]=0$ in $H_1(\mathbb{T})$.

 We begin by identifying $S^3$ with the set of unit quaternions in $\R^4$ with standard
coordinates: 
\[
 (w,x,y,z)=w+x\,\mathbf{i}+y\,\mathbf{j}+z\,\mathbf{k},
\]
and choosing a specific parametrization of the
Borromean rings $L^\text{Bor}$ in $S^3$. That is, define the $L_1$ component of
 $L^\text{Bor}$ to be the great circle  in $S^3$ through $1$ and $\mathbf{k}$, parametrized as
\[
x(s)=\cos(s)+\sin(s)\,\mathbf{k},\quad
\bigl(\,x(s)^{-1}=\cos(s)-\sin(s)\,\mathbf{k}\,\bigr).
\]


\begin{figure}

\begin{picture}(300, 120)
  \put(0,0){\includegraphics[height=1.8in,width=\textwidth]{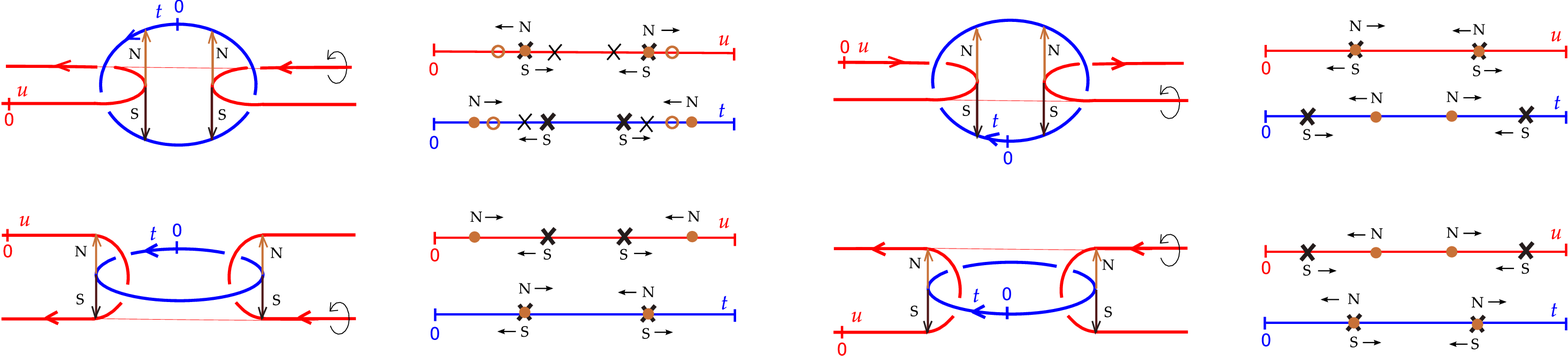}}
  \put(250,50){$s=\frac{3\pi}{2}$}
  \put(250,120){$s=\pi$}
  \put(0,60){$s=\frac{\pi}{2}$}
  \put(0,120){$s=0$}
\end{picture}
\caption{Four positions corresponding to angles: $s=0,\frac{\pi}{2}, \pi, \frac{3\pi}{2}$, small arrows next to $N$ and $S$ indicate a motion as $s \nearrow$.} \label{fig:4-pos}

\end{figure}

\no Observe that $\text{pr}(x(s))$ parameterizes the $z$-axis in $\R^3$. 
Figure \ref{fig:conf1} shows how to define the second and the third component $\{L_2, L_3\}$ of the Borromean rings $L^\text{Borr}$ in $\R^3$ considered as an image of $S^3-\{1\}$ under the stereographic projection $\text{pr}:S^3\subset \R^4\longrightarrow \R^3$ from $1\in S^3$. 
$L_2$ will bound the annuli with a rounded wedge removed, i.e. an arc $A_1$ of the circle
of radius $1$. The arc  $A_2$  belongs to the circle of radius $r_\epsilon=(1+\epsilon)$ in the $(x,y)$-plane. 
 The component $L_3$ is chosen to be a vertical ellipse
linking with $L_2$.

\begin{figure}
\begin{picture}(390, 250)(-5,0)
  \put(0,0){\includegraphics[width=\textwidth]{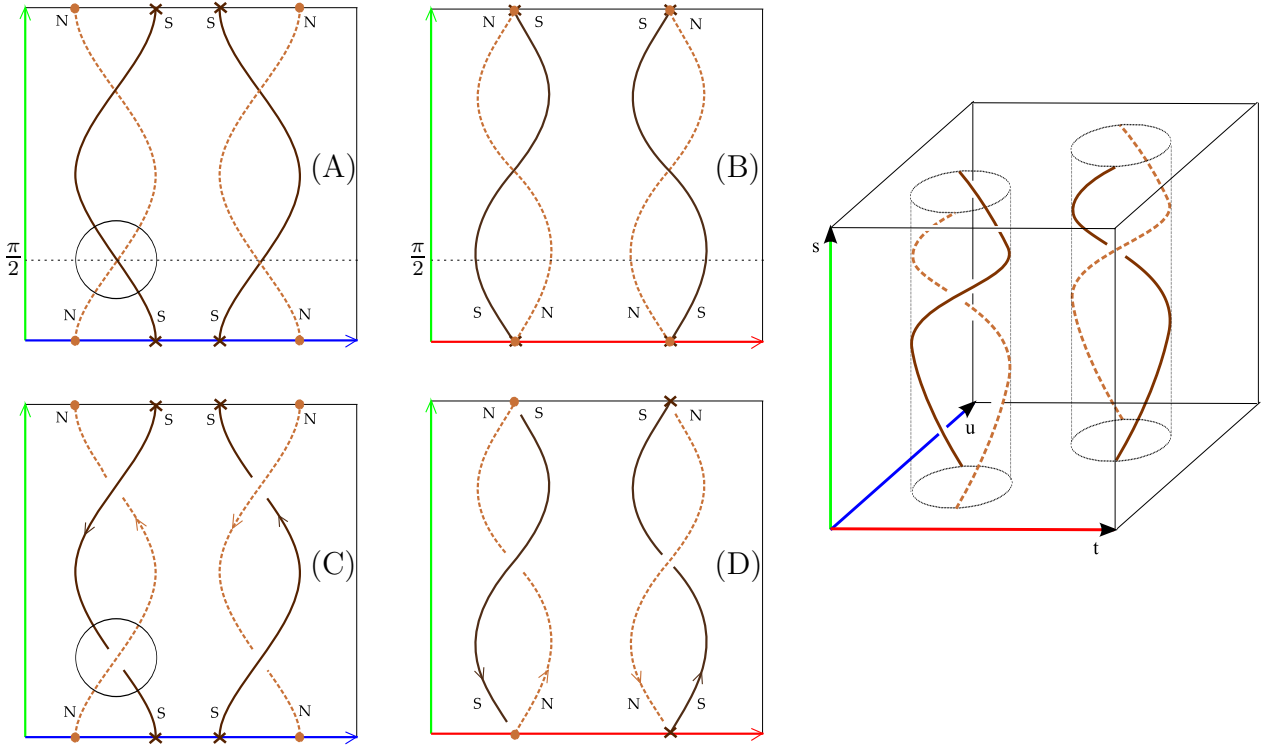}}
  \put(-6,182){$\frac{\pi}{2}$}
  \put(147,182){$\frac{\pi}{2}$}
  \put(110,215){(A)}
  \put(263,215){(B)}
  \put(110,65){(C)}
  \put(263,65){(D)}
\end{picture}
\caption{Projection of $l_{S,N}$ on the $su$-face and $st$-face of $\mathbb{T}$. The strands $l_S$ (solid line) and $l_N$ (dashed line) are oppositely oriented since $l_S$ and $l_N$ are null homologous in $\mathbb{T}$.}
\label{fig:framed-link}
\end{figure}

Next we focus on the Formula \eqref{eq:F_L}, observe that multiplication by $x(s)^{-1}$ has an effect of a rotation by angle $s$
in $(w,z)$-plane and $(x,y)$-plane of $\R^4$, which can be directly calculated:
\[
\begin{split}
x(s)^{-1}\cdot(w,x,y,z) & = \bigl(\cos(s) w+\sin(s) z, \cos(s) x +\sin(s) y,\\
&  \qquad \cos(s) y-\sin(s) x, \cos(s) z-\sin(s) w\bigr)\ .
\end{split}
\]

\no The flow defined by this $S^1$-action is tangent to the great circles of $S^3$, thus 
the projected flow on $\R^3$, via the stereographic projection $\text{pr}$,
presents the standard picture of the Hopf fibration. Let us call an invariant Hopf torus an $r$-torus, if and only if, it contains a circle of radius $r$ in the $(x,y)$-plane. Without loss of generality we assume that $L_2$ on Figure \ref{fig:conf1} belongs to the $r_{\epsilon/2}$-Hopf torus.  
 Every point on a $r$-torus traces a
$(1,1)$-curve under the $S^1$-action. For sufficiently small $\epsilon$, this motion can be regarded 
as a composition of the rotation by angle $s$ in both the direction of the meridian
and the longitude of a $r$-torus.  Therefore, for different values of $s$ the $S^1$-action ``rotates'' the components $L_2$ and $L_3$, by sliding 
along the Hopf tori by angle $s$ in the meridian and the longitudinal
direction. We denote resulting link components by 
\[
 L^s_2(t)=\text{pr}(x^{-1}(s)\cdot y(t)),\qquad \text{and}\qquad  L^s_3(u)=\text{pr}(x^{-1}(s)\cdot z(u)).
\]
The  unit circle on $(x,y)$-plane is left invariant under this action and therefore can be considered as the ``axis of the rotation''. This justifies the choice of the particular shape of $L^\text{Borr}$ pictured on Figure \ref{fig:conf1}.
\no Next, we seek to visualize the projection of $l_{S,N}=l_S\cup l_N$ on the $su$-face and $st$-face of the domain $\mathbb{T}$ of $F_{L^\text{Bor}}$ parameterized by $(s,t,u)\in \mathbb{T}$, (it is convenient to  think about $\mathbb{T}$ as a cube in $(s,t,u)$-coordinates, see Figure \ref{fig:framed-link}). 
For example when $s=0$, ($x(0)=1$), a point
$(0,t_0,u_0)$ belongs to $l_N$, if and only if, the vector $\mathsf{v}_0=\text{pr}(y(t_0))-\text{pr}(z(u_0))$
points  in the direction of $N=(0,0,1)$, analogous condition holds for direction
$S=(0,0,-1)$, and $l_S$. In order to determine a diagram of $l_{S,N}$, we must keep track of the ``head'' and ``tail'' of the vector $\mathsf{v}_s=L^s_2(t)-L^s_3(u)$, for various values of $s$ and record values of $t$ and $u$
for which $\mathsf{v}_s$ points ``North'' and ``South''
(Figure \ref{fig:4-pos}). This reads as the following condition
\begin{equation}\label{eq:arrow-eq}
  (s,t,u)\in l_{S,N},\qquad\text{if and only if},\qquad \mathsf{v}_s \parallel S\ \text{or}\ N.
\end{equation}
Without loss of generality we assume that $L^{s}_2$ is parametrized by the unit $t$-interval, and $L^{s}_3$ is parametrized by the unit $u$-interval. 
The process of recording values of $u$ and $t$ such that \eqref{eq:arrow-eq} holds is self-explanatory and is shown on Figure
\ref{fig:conf1} for values $s=0, \frac{\pi}{2}, \pi,
\frac{3\pi}{2}$, which is sufficient to draw projections of $l_S$ and $l_N$ on $st$- and $tu$-faces of $\mathbb{T}$.  Collecting the information on Figure \ref{fig:4-pos}, we draw the projection of
$l_{S,N}$ on the $su$-face of $\mathbb{T}$ represented
by square (A) in Figure \ref{fig:framed-link}. Analogously, the projection of 
$l_{S,N}$ on the $st$-face of $\mathbb{T}$ is obtained  and pictured in square (B).
In order to obtain the diagram of
$l_{S,N}$ we resolve the double points of Diagram (A) into crossings. For example, let us resolve the 
``circled'' double point on (A), which occurs at $s=\frac{\pi}{2}$ in the \emph{left} two stands of $l_{S,N}$. It suffices to determine the value of the $t$-coordinate at this point.
Diagram (B) tells us that $l_S$ is below $l_N$, because $\mathbb{T}$ is oriented so that the $t$-axis points above the $su$-face (see Figure \ref{fig:framed-link}).
 Resolving the remaining crossings in
a similar fashion leads to a diagrams of $l_{S,N}$ presented in squares (C) and (D). 
 Clearly, the linking number of
$l_S$ and $l_N$ is equal to $\pm 2$ in Diagram (C), (as the intersection number of e.g. $l_S$ with the obvious annulus on Diagram (C)). This justifies \eqref{eq:hopf=2}, and ends the proof.
\end{proof}
 
\no Results of Theorem \ref{th:milnor-hopf} and Proposition \ref{th:hopf-degree-general} combined with Formula \eqref{eq:hopf-integral} allow us to express $\bar{\mu}_{123}(L)$ of a Borromean link $L$ as
 \begin{equation}\label{eq:mu-formula}
 \begin{split}
  \bar{\mu}_{123}(L) & =  \int_{\mathbb{T}} F^\ast_L\nu\wedge \alpha=\int_{\mathbb{T}} L^\ast\omega\wedge \alpha,\\
  & \text{for}   \quad F^\ast_L\nu=d\alpha,\quad \omega=H^\ast\nu\in \Omega^2(\Cnf_3(S^3)),
  \end{split}
 \end{equation}
 where $\nu$ is the area form on $S^2$, and $H:\Cnf_3(S^3)\longrightarrow S^2$ is the deformation retraction 
 (as e.g. in \eqref{eq:hermans-map}). Alternatively, we may view $\omega$ as a 2-form on $(S^3)^3$ which is singular along 
 the diagonals $\mathbf{\Delta}\subset (S^3)^3$, and the singularity is of order $O(r^2)$, where $r$ is a distance to $\mathbf{\Delta}$. Consequently, $\omega$ is integrable but not square integrable on $(S^3)^3$.

\begin{remark}\label{rem:anti-symmetry}
  Notice that the integral formula \eqref{eq:mu-formula} exhibits the following property of $\bar{\mu}_{123}$:
   \[
  \bar{\mu}_{123}(L_1,L_2,L_3)=\text{\rm sign}(\sigma)\bar{\mu}_{123}(L_{\sigma(1)},L_{\sigma(2)},L_{\sigma(3)}),\qquad \sigma\in\Sigma_3\ .
 \]
\end{remark}


\no \section{Invariants of volume preserving flows. Helicities.}\label{sec:invariants-flows}

Given finitely many  volume preserving vector fields $B_1$, $B_2$, $\ldots B_k\in \text{\rm SVect}(M)$  on $M=S^3$ or a homology 3-sphere one seeks quantities $\mathsf{I}(B_1,B_2,\ldots,B_k)$ invariant under the action of volumorphisms isotopic to the identity $g\in \text{\rm SDiff}_0(M)$, commonly known as \emph{helicities} or \emph{higher helicities}:
\begin{gather}\label{eq:invariants}
 \mathsf{I}(B_1,B_2,\ldots,B_k)=\mathsf{I}(g_\ast B_1,g_\ast B_2,\ldots,g_\ast B_k),\qquad\quad \text{for all }g\in \text{\rm SDiff}_0(M), 
\end{gather}
where $g_\ast$ is a push-forward by a diffeomorphism $g$. To distinguish the case of a single vector field $B$ (i.e. $B=B_1=\ldots=B_k$) we often refer to $\mathsf{I}(B)=\mathsf{I}(B,B,\ldots,B)$ as \emph{self helicity}. We elucidated in the introduction a fundamental example of such invariant is the ordinary helicity $\mathsf{H}(B_1,B_2)$ of a pair of vector fields. In the remaining part of this section we review 
well known formulations of the helicity, which will later help us to point out analogies to the 
proposed formulation of the 3rd order helicity. Let $\mathcal{T}=\mathcal{T}_1\cup \mathcal{T}_2$ represent two invariant subdomains (not necessarily disjoint) under flows of $B_1$ and $B_2$ in $S^3$ and let 
\[
 \mathcal{T}=\mathcal{T}_1\times\mathcal{T}_2\subset \Cnf_2(M)\subset M\times M.
\]
Recall that the formula for $\mathsf{H}(B_1,B_2)$, from \cite{Khesin05, Vogel03}, specialized to invariant subdomains $\mathcal{T}=\mathcal{T}_1\cup \mathcal{T}_2$ may be expressed as 
\begin{equation}\label{eq:helicity}
 \mathsf{H}_{12}(B_1,B_2)=\int_{\mathcal{T}_1\times\mathcal{T}_2} \omega\wedge\iota_{B_1}\mu\wedge\iota_{B_2}\mu,
\end{equation}
where $\omega$ is known as the linking form on $M\times M$. When $\mathcal{T}=M\times M$ this formula is equivalent to a more commonly known expression: $\mathsf{H}(B_1,B_2)=\int_M \iota_{B_1}\mu\wedge d^{-1}(\iota_{B_2}\mu)$ (because $\omega$ also represents the integral kernel of $d^{-1}$).
 Philosophically, $\mathsf{H}(B_1,B_2)$  can be derived (c.f. \cite{Arnold86}) from the linking number of a pair of closed curves, which is expressed by Arnold's Helicity Theorem. For orbits
 $\{\mathscr{O}^1(x),\mathscr{O}^2(y)\}$  of $B_1$ and $B_2$ through $x, y\in M$, we introduce the following notation for the long pieces of closed up orbits
 \begin{equation}\label{eq:orbits}
\begin{split}
 \mathscr{O}^{B_i}_T(x) & =  \{\Phi^i_t(x)\ |\ 0\leq t\leq T\}\subset \mathscr{O}^i(x),\qquad i=1,2\\
 \bar{\mathscr{O}}^i_T(x) & :=  \mathscr{O}^i_T(x)\cup \sigma(x,\Phi^i(x,T)).
\end{split}
\end{equation}
 where $\sigma(x,y)$ denotes a short path, \cite{Vogel03}, connecting $x$ and $y$ in $M$ (see Section \ref{sec:ergodic}). Paraphrasing \cite{Arnold86} we state (for the proof also see \cite{Vogel03}),
 
\begin{theorem}[Arnold's Helicity Theorem, \cite{Arnold86}]
 Given $B_1, B_2\in\text{\rm SVect}(M)$,  The following 
 limit exists almost everywhere on $M\times M$:
 \begin{gather}\label{eq:linking-function}
  \bar{m}_{B_1 B_2}(x,y)=\lim_{T\to \infty} \frac{1}{T^2}\, \text{lk}(\bar{\mathscr{O}}^1_T(x)\times\bar{\mathscr{O}}^2_T(y))\ .
 \end{gather}
Moreover, $\bar{m}_{B_1 B_2}$ is in $L^1(M\times M)$, and 
\begin{equation}\label{eq:helicity-asymp}
 \mathsf{H}_{12}(B_1,B_2)=\int_{M\times M} \bar{m}_{B_1 B_2}(x,y)\, \mu(x)\wedge\mu(y)\ .
\end{equation}
\end{theorem}
\no The function $\bar{m}_{B_1 B_2}$ represents an asymptotic linking number of orbits $\{\mathscr{O}^1(x),\mathscr{O}^2(y)\}$, and the identity \eqref{eq:helicity-asymp} tells us that 
the helicity $\mathsf{H}_{12}(B_1,B_2)$ is equal to the average asymptotic linking number. In coming paragraphs, we will demonstrate,  how this philosophy is applied to obtain the asymptotic $\bar{\mu}_{123}$-invariant for 
3-component links and the third order helicity. 


\section{Definition of ``\texorpdfstring{$\bar{\mu}_{123}$}{mu-123}-helicity'' on invariant unlinked handlebodies.}\label{sec:handlebody}

\no In this section we apply the formulation of the $\bar{\mu}_{123}$-invariant for the 3-component links in $S^3$, 
obtained in Section \ref{sec:mu-and-hopf}, to define the third order helicity of a volume preserving vector field $B$ on certain invariant sets $\mathcal{T}$ of $B$ in $S^3$. 
In the following paragraphs as a ``warm-up'' to a more general case treated in Section \ref{sec:ergodic},  we consider the case of three disjoint unlinked handlebodies: $\mathcal{T}=\mathcal{T}_1\cup 	\mathcal{T}_2\cup \mathcal{T}_3$ in $S^3$ each of genus $g(\mathcal{T}_i)$. Henceforth, we use ``unlinked'' to mean ``with pairwise unlinked connected components''. 
When $\mathcal{T}$ represents three unlinked tubes (also known as \emph{flux tubes} \cite{Khesin98}) in $\R^3$, the third order helicity has been developed by several authors \cite{Berger90,Mayer03,Laurence-Stredulinsky00b} via Massey product formula for the $\bar{\mu}_{123}$-invariant, we compare our approach to these known works in Section \ref{sec:massey}.


\begin{figure}[htbp]
	\begin{picture}(400,190)
		\put(10,30){\includegraphics[width=2in]{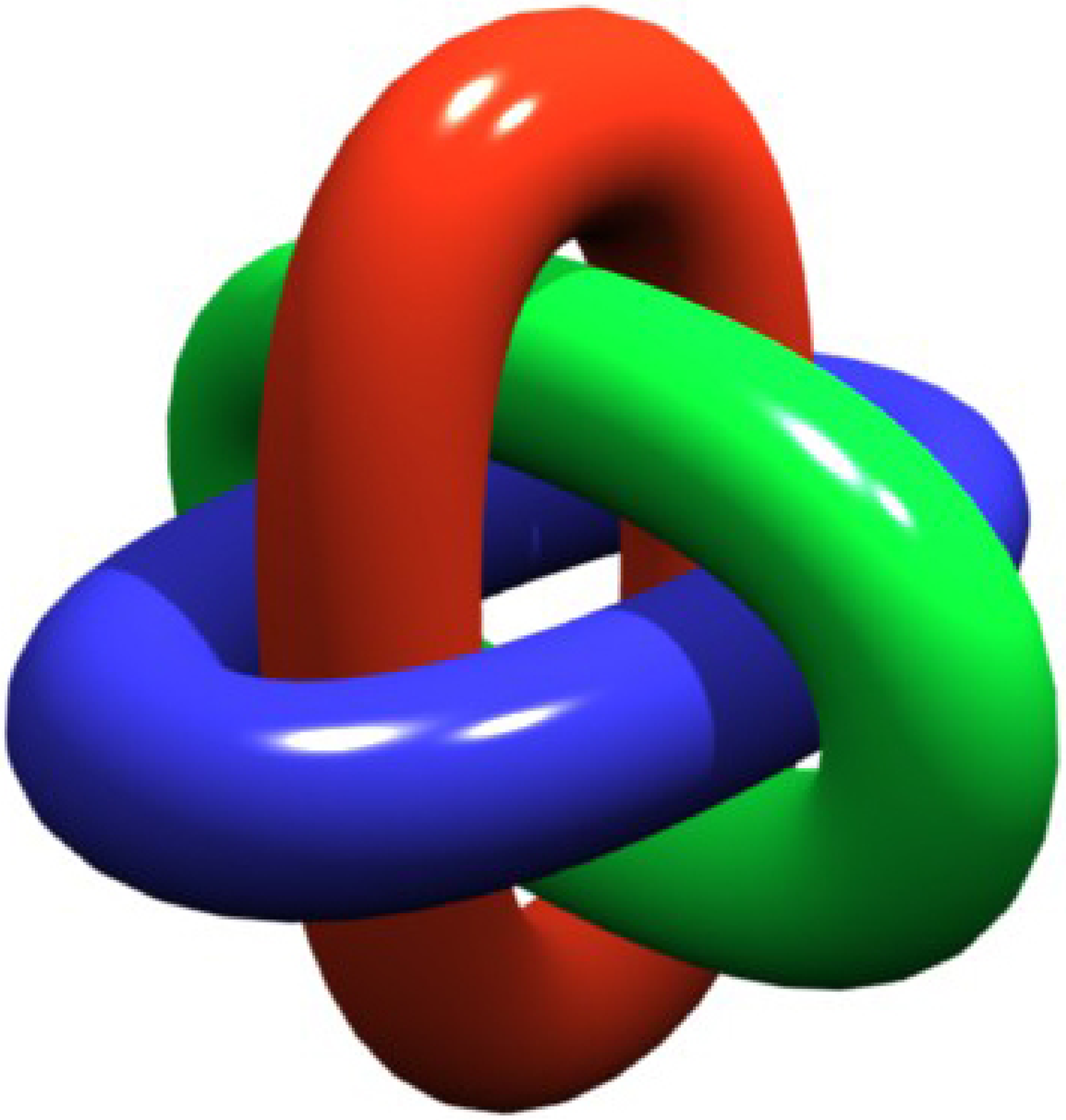}}
		\put(180,0){\includegraphics[width=4in]{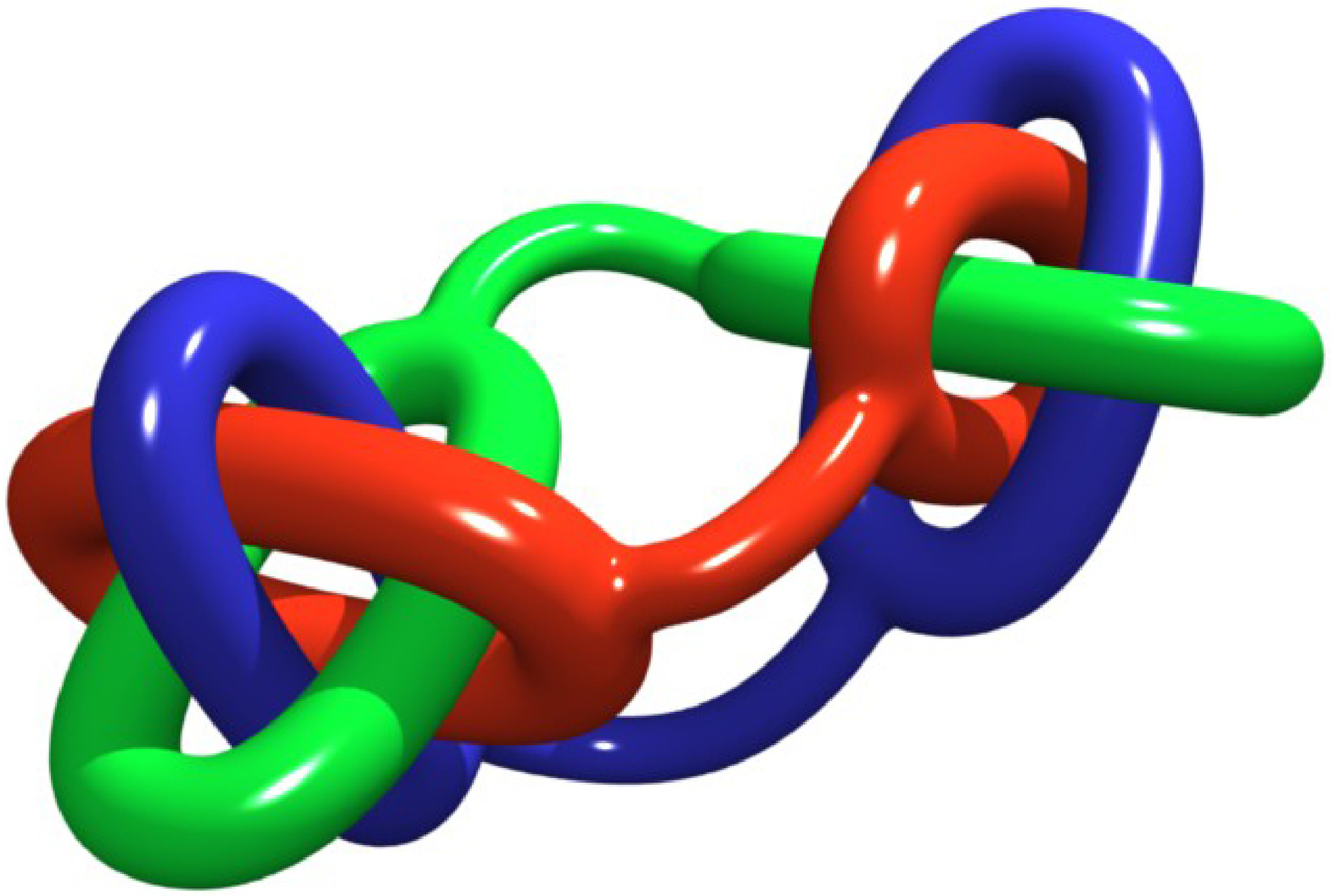}}
		\put(20,30){$\mathcal{T}^\text{Borr}$}
		\put(60,140){$\mathcal{T}_2$}
		\put(25,80){$\mathcal{T}_3$}
		\put(125,70){$\mathcal{T}_1$}
 \end{picture}
	\caption{The simplest unlinked handlebodies: flux tubes $\mathcal{T}^\text{Borr}$ modeled on the Borromean rings (left), and unlinked genus $2$ handlebodies (right).}
	\label{fig:invariant-sets}
\end{figure}

\no Assume $\mathcal{T}_i$s  have smooth boundary and $B$ to be tangent to $\partial \mathcal{T}_i$, we set
\[
 B_i:=B\bigl|_{\mathcal{T}_i},\qquad i=1,2,3,
\]
 and denote the flow of $B$ on $S^3$ by $\Phi$, and flows of restrictions $B_i$ by $\Phi^i$.  Clearly,
 such $\mathcal{T}$ is an invariant set of $B$.
 Given any domain $\mathcal{T}$ with three connected components $\{\mathcal{T}_i\}$ we may always associate a product domain in $\Cnf_3(S^3)$ as follows
\begin{equation}\label{eq:admissible-domain}
 \mathcal{T}:=\mathcal{T}_1\times \mathcal{T}_2\times \mathcal{T}_3\subset \Cnf_3(S^3)\subset S^3\times S^3\times S^3\ .
\end{equation}
Notice that $\mathcal{T}$ is a domain with corners in $\Cnf_3(S^3)$, and we use the same notation for the product of $\mathcal{T}_i$ as for the union in $S^3$. Wherever needed, we also assume that 
$(S^3)^3$ is equipped with a product Riemannian metric.
 Let a domain $\mathcal{T}$ defined in \eqref{eq:admissible-domain}, where $\mathcal{T}_i\cap \mathcal{T}_j=$\O, $i\neq j$ and each $\mathcal{T}_i$ is a handlebody in $S^3$ be called \emph{unlinked handlebody}, if and only if, the 2-form $\omega\in\Omega^2(\Cnf_3(S^3))$ defined in Equation \eqref{eq:mu-formula} is exact on $\mathcal{T}$, i.e. $\omega$ has a local \emph{potential} $\alpha_\omega\in \Omega^1(\mathcal{T})$: 
\begin{equation}\label{eq:omega-exact}
 \omega=d\alpha_\omega\ .
\end{equation}
Denote a volume preserving vector field $B$ and an unlinked handlebody $\mathcal{T}$ as a pair $(B;\mathcal{T})$.

\begin{remark}
\no Since $\omega$ is a dual cohomology class to the $S^2$ factor in $\Cnf_3(S^3)\cong S^3\times S^2$,  $\omega$ does not admit a global potential.
\end{remark}
\no Because each handlebody $\mathcal{T}_i$ has a homotopy type of a bouquet of circles, there is a natural choice of the basis for $H_1(\mathcal{T}_i)$ which consists of cycles $\{L^k_i\}_{k=1,\ldots,g(\partial\mathcal{T}_i)}$ corresponding to the circles. We have the following practical characterization of unlinked handlebodies:

\begin{lemma}\label{lem:unlinked-handlebody}
  $\mathcal{T}$ is an unlinked handlebody, if and only if, 
  \begin{equation}\label{eq:admissible-condition}
  \text{lk}(L^k_i,L^r_j)  = 0,\qquad \text{for all\ }\, i,j,\qquad i\neq j
 \end{equation}
\end{lemma}
\begin{proof}
The standard integral pairing, \cite{Bott82}, $H^2(\mathcal{T})\times H_2(\mathcal{T})\longrightarrow \R$ implies that a closed $k$-form is exact, if and only if, it evaluates to zero on all $k$-cycles of the domain. By the K$\ddot{\rm u}$nneth formula $H_2(\mathcal{T})=H_2(\mathcal{T}_1\times \mathcal{T}_2\times \mathcal{T}_3)$ is generated by $L^k_i\otimes L^r_j$, and by \eqref{eq:linking-numbers}:
\[
 \text{lk}(L^k_i,L^r_j) = \int_{L^k_i\times L^r_j} \omega\qquad  i\neq j.
\]
Therefore the condition \eqref{eq:admissible-condition} is necessary and sufficient for $\omega$ to be exact on $\mathcal{T}$.
\end{proof}

\no We define the $\bar{\mu}_{123}$-helicity of $(B;\mathcal{T})$  denoted by $\mathsf{H}_{123}(B;\mathcal{T})$ or $\mathsf{H}_{123}(B_1,B_2,B_3)$
as follows:
\begin{equation}\label{eq:mu_123-helicity}
\boxed{\mathsf{H}_{123}(B;\mathcal{T})=\mathsf{H}_{123}(B_1,B_2,B_3)\stackrel{\text{def.}}{=}\int_{\mathcal{T}} \bigl(\alpha_\omega\wedge d\alpha_\omega\bigr)\wedge \iota_{B_1}\mu_1\wedge \iota_{B_2}\mu_2\wedge \iota_{B_3}\mu_3, }
\end{equation}
where $\mu_i$ denotes the pull-back of the volume form $\mu$ on $S^3$ under the projection 
\begin{equation}\label{eq:projection-i}
\pi_i:S^3\times S^3\times S^3\longrightarrow S^3,\qquad \pi_i(x_1,x_2,x_3)=x_i,
\end{equation}
and $\iota_{B_i}$ is a contraction by a vector field $B_i$. Our notational convention is to denote by $\iota_{B_i}\mu_i$ both the forms on the base of $\pi_i$ and the pullbacks: $\pi^\ast_i\bigl(\iota_{B_i}\mu_i\bigr)$. Notice that $\mu=\mu_1\wedge\mu_2\wedge\mu_3$ is a volume form on the product: $S^3\times S^3\times S^3$.   
 There are obvious analogies between Formula \eqref{eq:mu_123-helicity} above, Formula \eqref{eq:helicity} for $\mathsf{H}_{12}(B;\mathcal{T})$, and the integral formula \eqref{eq:mu-formula} for the 
$\bar{\mu}_{123}$-invariant. The $3$-form:
 \[
 \gamma_\omega:=\alpha_\omega\wedge d\alpha_\omega=\alpha_\omega\wedge \omega,
 \]
 plays a role of the linking form as $\omega$ in Formula \eqref{eq:helicity}. 
 The main motivation behind definition \eqref{eq:mu_123-helicity} is the ergodic interpretation  of $\mathsf{H}_{123}(B;\mathcal{T})$ as an average asymptotic $\bar{\mu}_{123}$-invariant of orbits of $B$, which will become apparent in Section \ref{sec:ergodic}. Formula \eqref{eq:mu_123-helicity} can be also 
 regarded as the third order helicity of three distinct vector fields $B_i$, supported on the handlebodies $\mathcal{T}_i$. In Section \ref{sec:energy}, we indicate how to construct the potential $\alpha_\omega$ from the basic elliptic theory of differential forms.

\begin{theorem}[Helicity Invariance Theorem]\label{th:invariance-th}
 On every unlinked invariant handlebody $\mathcal{T}$ in $S^3$, $\mathsf{H}_{123}(B;\mathcal{T})$  is 
 \begin{itemize}
 \item[(i)] independent of a choice of the potential $\alpha_\omega$,
 \item[(ii)] invariant under the action of $\text{\rm SDiff}_0(S^3)$, i.e. for every $g\in \text{\rm SDiff}_0(S^3)$:
\begin{gather}\label{eq:invariance}
 \mathsf{H}_{123}(B;\mathcal{T})=\mathsf{H}_{123}(g_\ast B;g(\mathcal{T}))\ .
\end{gather}
 \end{itemize}
\end{theorem}
\begin{proof}
To prove $(i)$ observe for every other potential $\alpha'_\omega$ of $\omega$, the difference $\beta=\alpha_\omega-\alpha'_\omega$ is a closed $1$-form on $\mathcal{T}$ (since $\omega=d\alpha_\omega=d\alpha'_\omega$). Therefore,
\begin{eqnarray*}
 \mathsf{H}_{123}(B;\mathcal{T})-\mathsf{H}_{123}'(B;\mathcal{T}) & = & \int_{\mathcal{T}} (\beta\wedge\omega)\wedge\iota_{B_1}\mu_1\wedge \iota_{B_2}\mu_2\wedge \iota_{B_3}\mu_3\\
 & \stackrel{(1)}{=} & \int_{\mathcal{T}} d\bigl(\beta\wedge \alpha_\omega\wedge\iota_{B_1}\mu_1\wedge \iota_{B_2}\mu_2\wedge \iota_{B_3}\mu_3\bigr)\\
 & = & \int_{\partial\mathcal{T}} \beta\wedge \alpha_\omega\wedge\iota_{B_1}\mu_1\wedge \iota_{B_2}\mu_2\wedge \iota_{B_3}\mu_3 \stackrel{(2)}{=} 0,
\end{eqnarray*}
where in (1) we applied $d(\iota_{B_i}\mu_i)=0$ (since $B_i$'s are divergence free), and 
in (2):
\[
 \iota_{B_i}\mu_i\bigl|_{\partial\mathcal{T}_i}=0,
\]
(because each vector field $B_i$ is tangent to the boundary $\partial \mathcal{T}_i$), where
\[
 \partial \mathcal{T}=\bigl(\partial \mathcal{T}_1\times \mathcal{T}_2\times \mathcal{T}_3\bigr)\cup \bigl(\mathcal{T}_1\times \partial\mathcal{T}_2\times \mathcal{T}_3\bigr)\cup \bigl(\mathcal{T}_1\times \mathcal{T}_2\times \partial\mathcal{T}_3\bigr).
\]

\no The proof of $(ii)$ is in the style of 
\cite{Berger90, Mayer03}, but adapted to our setting. For any given $g\in \text{\rm SDiff}(S^3)$, by definition, there exists a path $t\longrightarrow g(t)\in \text{\rm SDiff}_0(S^3)$,
such that
\[
 g(0)=\text{id}_{S^3},\qquad g(1)=g\ .
\]
 Denote by $V$ the divergence free vector field on $S^3$, given by $V(x)=\frac{d}{dt} g(t,x)|_{t=0}$, i.e. $g(t)$ is a flow of $V$, and push-forward fields $B_i$ by 
\[
 B^t_i:=g(t)_\ast B_i\ .
\]
It is well known (see Appendix \ref{apx:A}, or \cite[p. 224]{Freedman91-2})
 that 2-forms: $\iota_{B^t_i}\mu$ are frozen in the flow of $V$, i.e.
\begin{equation}\label{eq:forms-frozen-in}
 \frac{d}{dt}\bigl(g(t)^\ast\iota_{B^t_i}\mu\bigr) =(\partial_t+\mathcal{L}_V) \iota_{B^t_i}\mu = 0.
\end{equation} 
We also have a path $\hat{g}(t)=(g(t),g(t),g(t))$ in $\text{\rm SDiff}_0(S^3\times S^3\times S^3)$, which analogously leads to the vector field $\hat{V}=(V,V,V)$. (Recall that a tangent bundle $T(S^3)^3$ has a natural product structure). Equation \eqref{eq:forms-frozen-in} implies
\begin{equation}\label{eq:forms-frozen-in2}
(\partial_t+\mathcal{L}_{\hat{V}}) \bigl(\pi^\ast_i\iota_{B^t_i}\mu\bigr)=(\partial_t+\mathcal{L}_{\hat{V}}) \iota_{B^t_i}\mu_i= 0\ .
\end{equation} 
(In the second equation we merely revoke our notational conventions: $\iota_{B^t_i}\mu_i\equiv\pi^\ast_i\bigl(\iota_{g(t)^\ast B_i}\mu\bigr)$).

\no Let $\mathcal{T}(t)=\hat{g}(t)(\mathcal{T}(0))\subset \Cnf_3(S^3)$, we must show $\frac{d}{dt} 
\mathsf{H}_{123}(B^t_1,B^t_2,B^t_3)=0$. Notice that 
for small enough $\epsilon$ and $t\in (t_0-\epsilon,t_0+\epsilon)$ we can assume, by $(i)$, that $\alpha_\omega$ is a time independent potential obtained from slightly bigger domain $\widetilde{\mathcal{T}}$ which deformation retracts on $\mathcal{T}(t_0)$, and satisfies 
\[
  \mathcal{T}(t)\subset\widetilde{\mathcal{T}},\qquad \text{\rm for}\quad t\in (t_0-\epsilon,t_0+\epsilon)\ .
\]
 Without loss of generality set 
$t_0=0$, and $\hat{g}(0)=\text{id}_{(S^3)^3}$, at $t_0$ we calculate:
\begin{eqnarray*}
 \frac{d}{dt} 
 \Bigl(\mathsf{H}_{123}(g(t)^\ast B_1,g(t)^\ast B_2,g(t)^\ast B_3)\Bigr) & = &
 \frac{d}{dt} \int_{\mathcal{T}(t)} \alpha_\omega\wedge d\alpha_\omega\wedge \iota_{B^t_1}\mu_1\wedge \iota_{B^t_2}\mu_2\wedge \iota_{B^t_3}\mu_3\\
 & = & \int_{\mathcal{T}(0)} \frac{d}{dt} \hat{g}(t)^\ast\bigl(\alpha_\omega\wedge d\alpha_\omega\wedge \iota_{B^t_1}\mu_1\wedge \iota_{B^t_2}\mu_2\wedge \iota_{B^t_3}\mu_3\bigr)\\
 & = & \int_{\mathcal{T}(0)} \bigl(\mathcal{L}_{\partial_t+\hat{V}}(\alpha_\omega\wedge d\alpha_\omega)\bigr)\wedge \iota_{B^t_1}\mu_1\wedge \iota_{B^t_2}\mu_2\wedge \iota_{B^t_3}\mu_3,
\end{eqnarray*}
where in the last identity we applied \eqref{eq:forms-frozen-in2} and the product rule for the Lie derivative.
Now because $\omega\wedge\alpha_\omega$ is time independent (for $t\in (t_0-\epsilon,t_0+\epsilon)$), Cartan magic formula yields
\begin{eqnarray*}
 \mathcal{L}_{\partial_t+\hat{V}}(\alpha_\omega\wedge d\alpha_\omega) & = & \mathcal{L}_{\hat{V}}(\alpha_\omega\wedge d\alpha_\omega) = \iota_{\hat{V}} d(\alpha_\omega\wedge d\alpha_\omega)+ d(\iota_{\hat{V}} (\alpha_\omega\wedge d\alpha_\omega))\\
 & =  & d(\iota_{\hat{V}} (\alpha_\omega\wedge d\alpha_\omega)),
\end{eqnarray*}
where $d(\alpha_\omega\wedge d\alpha_\omega)=\omega\wedge\omega=0$. 
 Since $B^t_i$ are tangent to the boundary of $\mathcal{T}_i(t)$, the same argument as in the proof of $(i)$ shows that the right hand side of the previous equation vanishes.
\end{proof}

\begin{remark}
 Notice that the above argument indicates that if we replace $\alpha_{\omega}\wedge d\alpha_\omega$ by virtually any closed 3-form $\eta$ on $\Cnf_3(S^3)$ (or $(S^3)^3$) we obtain some invariant under frozen-in-field deformations. If $\eta$ is exact we obtain trivial invariants, therefore the only sensible candidates here are cohomology classes of $\Cnf_3(S^3)\cong S^3\times S^2$. In dimension $3$ it leaves us with a dual to the $S^3$ factor in $\Cnf_3(S^3)$. Based on the considerations in Section \ref{sec:mu-and-hopf}
  one may argue that an invariant obtained this way is trivial. Indeed, the cohomology class $\eta$ evaluated on any 3-torus obtained from a 3-component link in $S^3$ via the map $L$ in \eqref{eq:F}
  is zero. Therefore, one could apply the ergodic approach of Section \ref{sec:ergodic} to show that $\int \eta\wedge\iota_{B}\mu_1\wedge \iota_{B}\mu_2\wedge \iota_{B}\mu_3$ 
 defines a trivial invariant.
  The crucial obstacle in extending the formula in \eqref{eq:mu_123-helicity} to encompass the whole $(S^3)^3$ is the fact that the potential $\alpha_\omega$ cannot be globally defined on $\Cnf_3(S^3)$.
\end{remark}

\section{The ergodic interpretation of  \texorpdfstring{$\mathsf{H}_{123}(B;\mathcal{T})$}{H-123}}\label{sec:ergodic}

The following statement is often seen in literature \cite{Cantarella-DeTurck-Gluck01, Cantarella-DeTurck-Gluck-Teytel00}: 

\begin{quote}
\emph{Helicity measures the extent to which 
vector fields twist and coil around each other.}
\end{quote}

\no A beauty of Arnold's ergodic approach 
to the helicity $\mathsf{H}_{12}(B)$ is that it makes this statement precise, by interpreting 
$\mathsf{H}_{12}(B)$ as an average asymptotic linking number of orbits of $B$. But, it also has a
practical application as it allows us to extend our approach to certain invariant sets of $B$. In this section we apply 
this philosophy to our newly defined invariant $\mathsf{H}_{123}(B;\mathcal{T})$, and interpret it
as the average asymptotic $\bar{\mu}_{123}$-invariant of orbits of $B$ in $\mathcal{T}$. 
Moreover, this ergodic interpretation leads us to an alternative, more intuitive proof of Helicity Invariance Theorem \ref{th:invariance-th}.

We begin by observing that given a volume preserving vector field $B$ on $M$ and its flow $\Phi_t$, 
we may regard $B$ as three vector fields on $(M)^3$. Thus, $(\Phi,\Phi,\Phi)$ induces a natural $\R^3$ action defined as follows:
\begin{equation}\label{eq:R3-action}
 \mathbf{\Phi}:\R^3\times (M)^3\longrightarrow (M)^3,\qquad \bigl((s,t,u),x,y,z\bigr)\stackrel{\mathbf{\Phi}}{\longrightarrow} (\Phi(s,x),\Phi(t,y), \Phi(u,z)).
\end{equation}
Observe that $\mathbf{\Phi}$ is a volume preserving action on $(M)^3$.  Our analysis is rooted in techniques developed in \cite{Arnold86, Laurence-Avellaneda93, Laurence-Stredulinsky00b, Vogel03}, the main tool is the following

\begin{theorem}[Multi-parameter Ergodic Theorem, \cite{Becker81}]\label{th:multi-ergodic}
 For any real valued $L^1$-function $F$ on $(M)^3$, the time
averages under the action in \eqref{eq:R3-action}:
\begin{gather*}
\bar{F}(x,y,z)=\lim_{T\to \infty}\frac{1}{T^3}\int^T_0\int^T_0\int^T_0
 F(\Phi(x,s),\Phi(y,t),\Phi(z,u))\,d s\,d t\,d u
\end{gather*}
converge almost everywhere. In addition, the limit function
$\bar{F}$ satisfies
\begin{itemize}
 \item[(i)] $\|\bar{F}\|_{L^1((M)^3)}\leq \|F\|_{L^1((M)^3)}$,
 \item[(ii)] $\bar{F}$ is invariant under the $\mathbf{\Phi}$-action,
 \item[(iii)] if $(M)^3$ is of finite volume then
\begin{gather}\label{eq:ergodic-integrals-equal}
 \int_{(M)^3}\bar{F} =\int_{(M)^3} F\ .
\end{gather}
\end{itemize}
\end{theorem}

\begin{definition}

Define \emph{invariant unlinked domain} $\mathcal{T}$ of $B$ as an arbitrary $\mathbf{\Phi}$-invariant set, with topological closure $\overline{\mathcal{T}}$ which belongs to a larger product of open sets $\widetilde{\mathcal{T}}=\widetilde{\mathcal{T}}_1\times \widetilde{\mathcal{T}}_2\times \widetilde{\mathcal{T}}_3$
in $\Cnf_3(S^3)$, satisfying the following 
\begin{itemize}
\item[(A)] $\widetilde{\mathcal{T}}$ admits a \emph{short path system} $\mathcal{S}$,
\item[(B)] Equation \eqref{eq:omega-exact} holds on $\widetilde{\mathcal{T}}$,
\end{itemize}
\no where by a system of short paths on $\widetilde{\mathcal{T}}$, \cite{Vogel03, Khesin98}, we understand 
a collection of  curves $\mathcal{S}=\{\sigma_i(x,y)\}$ on each open set $\widetilde{\mathcal{T}}_i$ such that 
 
\begin{itemize}
\item[(a)] for every pair of points $x,y \in \widetilde{\mathcal{T}}_i$ there is a connecting curve $\sigma_i(x,y):I\mapsto \widetilde{\mathcal{T}}_i$ in $\mathcal{S}$, $\sigma_i(0)=x$ and $\sigma_i(1)=y$,
\item[(b)] the lengths of paths in $\mathcal{S}$ are uniformly 
bounded above by a common constant. 
\end{itemize}
\end{definition}

\no Topologically, every $\mathbf{\Phi}$-invariant 
set is a union of products of orbits of $B$ in $(S^3)^3$. It is often convenient to think of the 
orbits $\mathbf{\Phi}$-action as a foliation of $(S^3)^3$. Then $\mathbf{\Phi}$-invariant 
sets are just union of leaves of this foliation.  A fundamental example of an 
invariant unlinked domain is the case of $\mathbf{\Phi}$-invariant 
set $\mathcal{T}$ contained in the product $\widetilde{\mathcal{T}}=\widetilde{\mathcal{T}}_1\times \widetilde{\mathcal{T}}_2\times \widetilde{\mathcal{T}}_3$ of disjoint open unlinked handlebodies $\widetilde{\mathcal{T}}_i$. Note that in this case we do not require $B$ to be tangent to $\partial\widetilde{\mathcal{T}}_i$, and $\widetilde{\mathcal{T}}$ always admits a short path system as we describe in the following

\begin{remark}
In \cite{Vogel03}, Vogel shows that on a closed manifold $M$, geodesics always provide a short path system. When $\mathcal{T}$ is contained in the product of unlinked handlebodies $\widetilde{\mathcal{T}}$
we may easily construct such system on $\widetilde{\mathcal{T}}$ as follows
 Because  $\mathcal{\mathcal{T}}_i$ are proper subsets of $S^3$ we generally do not want to use ambient geodesics from $S^3$ as they may not lie entirely 
in $\mathcal{T}_i$. To obtain $\mathcal{S}$ one puts an artificial Riemannian metric 
on each $\widetilde{\mathcal{T}}_i$ which makes $\partial\widetilde{\mathcal{T}}_i$ totally geodesic, and choose $\mathcal{S}$ to be geodesics on such Riemannian manifold. Observe that applying a diffeomorphism $g\in \text{Diff}(\widetilde{\mathcal{T}})$ to $\mathcal{S}$ 
results in the system $g\mathcal{S}$ on $g(\widetilde{\mathcal{T}})$.
\end{remark}

\no The following result is an analog of Arnold's Helicity Theorem in our setting,

\begin{theorem}[Ergodic interpretation of $\mathsf{H}_{123}(B;\mathcal{T})$]\label{th:ergodic-mu}
 Given $(B;\mathcal{T})$, the following limit (asymptotic $\bar{\mu}_{123}$-invariant of orbits) exists for  almost 
 all $(x,y,z)\in \mathcal{T}$:
\begin{equation}\label{eq:bar-m_B}
 \bar{m}_B(x,y,z)=\lim_{T\to\infty} \frac{1}{T^3}\bar{\mu}_{123}\bigl(\bar{\mathscr{O}}^{B_1}_T(x),\bar{\mathscr{O}}^{B_2}_T(y),\bar{\mathscr{O}}^{B_3}_T(y)\bigr).
\end{equation}
Moreover, 
\begin{equation}\label{eq:H_123-ergodic}
 \mathsf{H}_{123}(B;\mathcal{T})=\int_{\mathcal{T}} \bar{m}_B(x,y,z)\,\mu(x)\wedge\mu(y)\wedge\mu(z)\ .
\end{equation}
\end{theorem}
\begin{proof}
 The proof is similar to the one in e.g. \cite{Vogel03}. Before we start, we must point out
 the following identity (valid for any 3-form $\beta$ on $M\times M\times M$ and vector fields 
 $B_1,B_2,B_3$ on $M$)
 
 \begin{equation}\label{eq:iota-volume}
 \begin{split}
(\iota_{B_3}\iota_{B_2}\iota_{B_1}\beta)\wedge\mu_1\wedge\mu_2\wedge\mu_3 & =
\beta(B_1,B_2,B_3)\,\mu_1\wedge\mu_2\wedge\mu_3\\
 & =\beta\wedge\iota_{B_1}\mu_1\wedge\iota_{B_2}\mu_2\wedge\iota_{B_3}\mu_3\ .
 \end{split}
 \end{equation}
  \no The first equation follows from the definition, the second one is a consequence of the fact that $\iota_B$ is an antiderivation i.e.
 \begin{equation}\label{eq:iota-antiderivation}
  \iota_B(\alpha\wedge\beta)=(\iota_B\alpha)\wedge\beta+(-1)^{|\alpha|}\alpha\wedge(\iota_B\beta),
 \end{equation}
 and $\iota_{B_i}\mu_j=0$, for $i\neq j$ (see Appendix \ref{apx:A}). As a result,
 \begin{eqnarray*}
\notag  \mathsf{H}_{123}(B_1,B_2,B_3) & = & \int_{\mathcal{T}} \bigl(\alpha_\omega\wedge d\alpha_\omega\bigr)\wedge \iota_{B_1}\mu_1\wedge \iota_{B_2}\mu_2\wedge \iota_{B_3}\mu_3\\
& = & \int_\mathcal{T} \bigl(\iota_{B_3}\iota_{B_2}\iota_{B_1}\bigl(\alpha_\omega\wedge d\alpha_\omega\bigr)\bigr)\,\mu_1\wedge\mu_2\wedge\mu_3\\
\label{eq:m-function}  & = & \int_\mathcal{T} m_{B_1,B_2,B_3}(x,y,z)\,\mu_1\wedge \mu_2\wedge \mu_3,
\end{eqnarray*}
where $m_{B_1,B_2,B_3}:=\alpha_\omega\wedge d\alpha_\omega(B_1,B_2,B_3)$. For convenience, let us set (see \eqref{eq:orbits}):
\[
 \bar{\mathscr{O}}(x,y,z;T):=\bar{\mathscr{O}}^{B_1}_T(x)\times \bar{\mathscr{O}}^{B_2}_T(y)\times\bar{\mathscr{O}}^{B_3}_T(y)\ .
\]
\no Observe that if the orbit $\bar{\mathscr{O}}(x,y,z;T)$ is nondegenerate, it represents a Borromean link and we may apply Formula \eqref{eq:mu-formula} to get
\begin{eqnarray*}
\bar{\mu}_{123}\bigl(\bar{\mathscr{O}}^{B_1}_T(x),\bar{\mathscr{O}}^{B_2}_T(y),\bar{\mathscr{O}}^{B_3}_T(y)\big)
& = & \int_{\bar{\mathscr{O}}(x,y,z;T)}
\alpha_\omega\wedge d\alpha_\omega\\
& = & \int_{\mathscr{O}(x,y,z;T)}
\alpha_\omega\wedge d\alpha_\omega+(I),
\end{eqnarray*}
where the term $(I)$ involves integrals over short paths in $\mathcal{S}$, (see Appendix \ref{apx:B}).
For degenerate orbits (such as fixed points etc.), the above formula still makes sense because $\bar{\mu}_{123}$ is a homotopy invariant of the associated map as proven in Section \ref{sec:mu-and-hopf}.
For every $(x,y,z)$ : 
\[
 D\mathbf{\Phi}_{(x,y,z)}[\partial_i]=B_i,\quad i=1,2,3,
\] 
where $\partial_1:=\partial_s$,$\partial_2:=\partial_t$,$\partial_3:=\partial_u$, thus we obtain
\begin{eqnarray*}
 \mathbf{\Phi}^\ast_{(x,y,z)}(\alpha_\omega\wedge d\alpha_\omega) & = & (\alpha_\omega\wedge d\alpha_\omega(B_1,B_2,B_3)(\Phi^1(x,s),\Phi^2(y,t),\Phi^3(z,u))\,ds\wedge dt\wedge du\\
 & = & m_{B_1,B_2,B_3}(\Phi^1(x,s),\Phi^2(y,t),\Phi^3(z,u))\,ds\wedge dt\wedge du\ .
\end{eqnarray*}
Therefore, 
\[
\int_{\mathscr{O}(x,y,z;T)}
\alpha_\omega\wedge d\alpha_\omega=\int^T_0\int^T_0\int^T_0 m_{B_1,B_2,B_3}(\Phi^1(x,s),\Phi^2(y,t),\Phi^3(z,u))\,ds\wedge dt\wedge du\ .
\]
Function $m_{B_1,B_2,B_3}$ is smooth bounded on $\mathcal{T}$ and hence $L^1$. Because short paths do not contribute to the time average (see Appendix \ref{apx:B}):
\[
 \lim_{T\to\infty} \frac{1}{T^3}(I)=0.
\]
Theorem \ref{th:multi-ergodic} applied to the function $m_{B_1,B_2,B_3}$ yields almost everywhere existence of the limit \eqref{eq:bar-m_B}. Hence we obtain the invariant $L^1$-function $\bar{m}_B:=\bar{m}_{B_1,B_2,B_3}$ on $\mathcal{T}$. The identity \eqref{eq:H_123-ergodic} 
follows from $(iii)$ of Theorem \ref{th:multi-ergodic}.
\end{proof}

\begin{theorem}[Helicity Invariance Theorem-ergodic version]\label{th:invariance-th-ergodic}
 On every unlinked domain $\mathcal{T}$ in $S^3$, (i) and (ii) of Theorem 
 \ref{th:invariance-th} hold for $\mathsf{H}_{123}(B;\mathcal{T})$.
 \end{theorem}

\begin{proof}
The proof of $(i)$ immediately follows from independence of the limit \eqref{eq:bar-m_B}
of the choice of the potential $\alpha_\omega$.
For the proof of  $(ii)$ we must show the following
\[
 \int_\mathcal{T} m_{B_1,B_2,B_3}\,\mu_1\wedge \mu_2\wedge \mu_3=\int_{g(\mathcal{T})} m_{g_\ast B_1,g_\ast B_2, g_\ast B_3}\,\mu_1\wedge \mu_2\wedge \mu_3\ .
\]
Theorem \ref{th:multi-ergodic} tells us that $m_{B_1,B_2,B_3}$ and $m_{g_\ast B_1,g_\ast B_2, g_\ast B_3}$ admit $L^1$-averages $\bar{m}_{B_1,B_2,B_3}$ and $\bar{m}_{g_\ast B_1,g_\ast B_2, g_\ast B_3}$, under actions of $B_i$ and $g_\ast B_i$ respectively. It suffices to show the following identity
\begin{equation}\label{eq:m=gm}
 \bar{m}_{B_1,B_2,B_3}(x,y,z)=\bar{m}_{g_\ast B_1,g_\ast B_2, g_\ast B_3}(g(x),g(y),g(z)),\qquad \text{a.e.}
\end{equation}
then \eqref{eq:invariance} is an immediate consequence of Equation \eqref{eq:ergodic-integrals-equal}, change of variables for integrals, and the fact that $g$ preserves volume (i.e. $\mu_i=g^\ast\mu_i$).
Borrowing notation from the previous theorem set
\begin{eqnarray*}
 g\bar{\mathscr{O}}(x,y,z;T) & := & g(\bar{\mathscr{O}}^{B_1}_T(x))\times g(\bar{\mathscr{O}}^{B_2}_T(y))\times g(\bar{\mathscr{O}}^{B_3}_T(z))\\
 & = & \bar{\mathscr{O}}^{g_\ast B_1}_T(g(x))\times \bar{\mathscr{O}}^{g_\ast B_2}_T(g(y))\times \bar{\mathscr{O}}^{g_\ast B_3}_T(g(z)),
\end{eqnarray*}
where the second identity is a consequence of the fact that the flow of $g_\ast B_i$ is 
obtained as a composition of $g$ and the flow of $B_i$.
Since $g$ is isotopic to the identity, $\bar{\mathscr{O}}(x,y,z;T)$ and $g\bar{\mathscr{O}}(x,y,z;T)$
 are homotopic as link maps (for nondegenerate orbits they are in fact isotopic  Borromean links in $S^3$)
 and by theorems of Section \ref{sec:mu-and-hopf}, we have
 \[
  \bar{\mu}_{123}(\bar{\mathscr{O}}(x,y,z;T))=\bar{\mu}_{123}(g\bar{\mathscr{O}}(x,y,z;T)).
 \]
 As a result of the above identity and \eqref{eq:bar-m_B} we derive a.e.
\begin{eqnarray*}
  \bar{m}_{B_1,B_2,B_3}(x,y,z) & = & \lim_{T\to\infty}\frac{1}{T^3}\bar{\mu}_{123}(\bar{\mathscr{O}}(x,y,z;T))\\
   & = & \lim_{T\to\infty}\frac{1}{T^3}\bar{\mu}_{123}(g\bar{\mathscr{O}}(x,y,z;T))=\bar{m}_{g_\ast B_1,g_\ast B_2, g_\ast B_3}(g(x),g(y),g(z)).
\end{eqnarray*}
Notice that in the last equation we used the ``pushed forward'' short paths system: $g\mathcal{S}$.
Since the lengths of paths in $g\mathcal{S}$ are bounded as well, they do not contribute to the limit. This proves the identity \eqref{eq:m=gm}, and consequently \eqref{eq:invariance}.
\end{proof}

\no Notice that the above argument does not require the Stokes Theorem, and as such may lead to further generalizations. Clearly, for $\mathsf{H}_{123}(B;\mathcal{T})$ to be nontrivial $\mathcal{T}$ must be 
of nonzero measure.


\section{Flux formula for  \texorpdfstring{$\mathsf{H}_{123}(B;\mathcal{T})$}{H-123}.}

 The following formula is a well known property of the ordinary helicity $\mathsf{H}_{12}(B;\mathcal{T})$ of the flux tubes $\mathcal{T}$ modeled on a 2-component link $L=\{L_1,L_2\}$ (see e.g. \cite{Laurence-Avellaneda93, Cantarella00})

\begin{equation}\label{eq:helicity12-flux-tubes}
\mathsf{H}_{12}(B;\mathcal{T})=\mathsf{H}_{12}(B_1,B_2)=\text{lk}(L_1,L_2)\,\text{Flux}(B_1)\,\text{Flux}(B_2).
\end{equation}
\no Here we show an analogous property for $\mathsf{H}_{123}(B;\mathcal{T})$, when $\mathcal{T}$ is  an invariant unlinked  handlebody. Recall that $\{L^k_i\}_{k=1,\ldots,g(\partial\mathcal{T}_i)}$ denotes the basis of $H_1(\mathcal{T})$ defined in Lemma \ref{lem:unlinked-handlebody}.

\begin{proposition}
 $\mathsf{H}_{123}(B_1,B_2,B_3)$ on invariant unlinked handlebodies $\mathcal{T}$ satisfies the following formula
 \begin{equation}\label{eq:helicity-fluxes}
 \mathsf{H}_{123}(B_1,B_2,B_3)=\sum_{i, j, k} \bar{\mu}_{123}(L^i_1, L^j_2, L^k_3)\,\text{Flux}_{\Sigma_i}(B_1)\,\text{Flux}_{\Sigma_j}(B_2)\,\text{Flux}_{\Sigma_k}(B_3),
 \end{equation}
where $\{L^i_1\otimes L^j_2\otimes L^k_3\}$ is a basis of $H_3(\mathcal{T})$, $\text{Flux}(B_i)$ stands for the flux of $B_i$ through a cross sectional surface $\Sigma_k$ of $\mathcal{T}_i$, which represents the homology Poincar$\acute{\rm e}$ dual of $L^k_i$ in $H_2(\mathcal{T}_i,\partial\mathcal{T}_i)$.

\end{proposition}
\begin{proof}
 Recall that the flux $\text{Flux}_{\Sigma_k}(B_i)$ of a vector field $B_i$ though a cross-sectional surface
 $\Sigma_k$ in $\mathcal{T}_i$ is given by:
 \begin{equation}\label{eq:flux}
 \begin{split}
  \text{Flux}_{\Sigma_k}(B_i) & =  \int_{\Sigma_k} \iota_{B_i}\mu=\int_{\mathcal{T}_i} h_k\wedge\iota_{B_i}\mu\\
  & =  \int_{\mathcal{T}_i} \iota_{B_i}h_k\wedge\mu =\int_{\mathcal{T}_i} h_k(B_i),
  \end{split}
 \end{equation}
 where 1-forms $h_k$ represent cohomology Poincar$\acute{\rm e}$ duals of $\Sigma_k$, and we applied \eqref{eq:iota-antiderivation} in the third equation. For every closed curve $\gamma\subset \mathcal{T}$   $h_k$ satisfies 
 \[
  \int_{\gamma} h_k=\#(\gamma,\Sigma_k)=\text{deg}(\gamma, L_k),
 \]
 where $\text{deg}(\gamma, L^j_k)$ measures how many times $\gamma$ ``wraps around'' the cycle  $L^j_k$. 
 
 For simplicity, we first assume that $\mathcal{T}$ is modeled on a Borromean link
 $L=\{L_1,L_2,L_3\}$ (such as $\mathcal{T}^\text{Borr}$ on Figure \ref{fig:invariant-sets}). Then, $h=h_1\wedge h_2\wedge h_3\in \Omega^3(\mathcal{T})$  is a cohomology class dual to the cycle $L_1\otimes L_2\otimes L_3$ in $H_3(\mathcal{T})$. Define
 \[
  H:=\iota_{B_3}\iota_{B_2}\iota_{B_1} h=h(B_1,B_2,B_3),
 \]
 which is a smooth function, and let $\bar{H}$ be the time average of $H$ as in Theorem \ref{th:multi-ergodic}. It suffices to show 
 \begin{equation}\label{eq:H-equation}
   \bar{m}_{B_1,B_2,B_3}=\bar{\mu}_{123}(L)\,\bar{H},\qquad \text{a.e.}
 \end{equation}
  Analogously, as in the proof of Theorem \ref{th:invariance-th}, 
 Equation \eqref{eq:H-equation} immediately implies Formula \eqref{eq:helicity-fluxes}. 
 Assume the notation of Theorem \ref{th:invariance-th}, for given $T$ consider $\bar{\mathscr{O}}(x,y,z;T)$. 
Thanks to Theorem \ref{th:milnor-hopf} (see the first paragraph of the proof) we have 
 \begin{eqnarray*}
  \bar{\mu}_{123}(\bar{\mathscr{O}}(x,y,z;T)) & = & \bar{\mu}_{123}(L)\,\bigl(\text{deg}(\bar{\mathscr{O}}^{B_1}_T(x),L_1)\,\text{deg}(\bar{\mathscr{O}}^{B_2}_T(y),L_2)\,\text{deg}(\bar{\mathscr{O}}^{B_3}_T(x),L_3)\bigr)\\
  & = & \bar{\mu}_{123}(L)\int_{\bar{\mathscr{O}}(x,y,z;T)} h\ .
 \end{eqnarray*}
Therefore,
 \begin{eqnarray*}
  \bar{m}_{B_1,B_2, B_3}(x,y,z) & = & \lim_{T\to \infty} \frac{1}{T^3}\bar{\mu}_{123}(\bar{\mathscr{O}}(x,y,z;T))\\
  & = & \bar{\mu}_{123}(L)\lim_{T\to \infty} \frac{1}{T^3}\int_{\bar{\mathscr{O}}(x,y,z;T)} H\\
  & =  &\bar{\mu}_{123}(L) \bar{H}(x,y,z),
 \end{eqnarray*}
where the last equality is again the consequence of short paths not contributing to the limit.
From the product structure of $\mathcal{T}$ and $(iii)$ of Theorem \ref{th:multi-ergodic} we get
 \[
  \int_\mathcal{T} \bar{H} =\int_\mathcal{T} H=\bigl(\int_{\mathcal{T}_1} h_1(B_1) \bigr)\bigl(\int_{\mathcal{T}_2} h_2(B_2) \bigr)\bigl(\int_{\mathcal{T}_3} h_3(B_3) \bigr),
 \]
 which combined with Equation \eqref{eq:flux} for fluxes concludes the proof in the case of Borromean flux tubes.
The proof in the case of a general handlebody is analogous, once we show the following 
\begin{lemma}
 Let $\mathcal{O}=\{\mathcal{O}_1,\mathcal{O}_2, \mathcal{O}_3\}$ be a 3-component link in $S^3$ such that 
 $\mathcal{O}_i\subset \mathcal{T}_i$ for each $i$. Then $\mathcal{O}$ is Borromean and 
  \begin{equation}\label{eq:mu-in-T}
 \bar{\mu}_{123}(\mathcal{O})=\sum_{i, j, k} \bar{\mu}_{123}(L^i_1, L^j_2, L^k_3)\,\text{deg}(\mathcal{O}_1, L^i_1)\,\text{deg}(\mathcal{O}_2, L^j_2)\,\text{deg}(\mathcal{O}_3, L^k_1)\ .
 \end{equation}
\end{lemma}
\begin{proof}
 Thanks to the interpretation of the $\bar{\mu}_{123}$-invariant in Section \ref{sec:mu-and-hopf}, it is not only a link homotopy 
 invariant, but also a homotopy invariant of the associated map $F_{\mathcal{O}}$ defined in \eqref{eq:F}. Observe that  each component 
 $\mathcal{O}_i$ can be homotoped inside of its handlebody $\mathcal{T}_i$ to become a bouquet of circles $\widehat{\mathcal{O}}_i\cong S^1\vee S^1\vee\ldots\vee S^1$ so that each factor in $\widehat{\mathcal{O}}_i$ is a multiple of the cycle represented by $\{L^j_i\}$in $H_1(\mathcal{T}_i)$. As a result, we obtain the associated map
 $F_{\widehat{\mathcal{O}}}$, and
 \[
  \bar{\mu}_{123}(\mathcal{O})=\frac{1}{2}\mathscr{H}(F_{\widehat{\mathcal{O}}}).
 \]
 Interpreting $\mathscr{H}(F_{\widehat{\mathcal{O}}})$ as the intersection number and summing up
 intersection numbers we conclude \eqref{eq:mu-in-T}.
\end{proof}
\end{proof}

\no In the case of Borromean flux tubes, Formula \eqref{eq:helicity12-flux-tubes} reduces to
\begin{equation}\label{eq:helicity12-borro-flux-tubes}
\mathsf{H}_{123}(B_1,B_2,B_3)=\bar{\mu}_{123}(L)\,\text{Flux}_{\Sigma_1}(B_1)
\text{Flux}_{\Sigma_2}(B_2)\text{Flux}_{\Sigma_3}(B_3),
\end{equation}
where $\Sigma_i$ denotes homology Poincar$\acute{\rm e}$ duals to $L_i$ in $H_2(\mathcal{T}_i,\partial\mathcal{T}_i)$.

  Since the fluxes are invariant under frozen-in-field deformations, Formula \eqref{eq:helicity12-flux-tubes}
  is yet another proof of Theorem \ref{th:invariance-th} in the setting of invariant unlinked handlebodies.
  In \cite{Laurence-Stredulinsky00b} the authors develop the same formula for the Borromean flux tubes.
  This clearly must be the case, as we work with the same topological 
  invariants of links via a different approach. Additional advantage 
  of our formulation is that we do not have to separately 
  deal with null points of vector fields as in \cite{Laurence-Stredulinsky00b}. 


\section{Energy bound.}\label{sec:energy}

\no In this section we indicate how the quantity $\mathsf{H}_{123}(B;\mathcal{T})$ invariant under frozen-in-field deformations provides a lower bound for the  $L^2$-energy $E_2(B)$ of a volume preserving field $B$ on $M=S^3$. We restrict our considerations to the case of an invariant unlinked  handlebody  $\mathcal{T}$, defined in Section \ref{sec:handlebody}. For the notation used in this section see Appendix \ref{apx:C}.

Recall the definition
\begin{equation}\label{eq:energy}
 E_2(B)=\int_M |B|^2 = \|B\|^2_{L^2(M)}\ .
\end{equation}
The ordinary helicity $\mathsf{H}_{12}(B)$ provides a well known lower bound (see \cite[p. 123]{Khesin98}):
\[
  \frac{1}{\lambda_1}|\mathsf{H}_{12}(B)|\leq  E_2(B),
\]
where $\lambda_1$ is the first eigenvalue of the elliptic self adjoint operator $\ast\,d:\Omega^1(M)\mapsto \Omega^1(M)$
(known as the curl operator), $\ast$ denotes the Hodge star. Importance of such lower energy bounds,
 stems from an area of interest in the ideal magnetohydrodynamics, \cite{Priest84}, as this constrains the phenomenon of ``magnetic relaxation'', \cite{Freedman99}.
A need for higher helicities can be justified by the fact that one may easily produce examples of vector fields $B$ for which $\mathsf{H}_{12}(B)$ vanishes, but the energy of the field $B$ still cannot be relaxed. For example, consider a classical case of Borromean flux tubes $\mathcal{T}^\text{Borr}$ with $B$ smooth and vanishing outside the tubes. Furthermore, assume that orbits of $B$ are just ``parallel 
circles'' inside each tube. By taking $B=B_1+B_2+B_3$, bilinearity of $\mathsf{H}_{12}(\,\cdot\,,\,\cdot\,)$ on $\text{\rm SVect}(S^3)$ and disjoint supports of $B_i$  yield
\[
 \mathsf{H}_{12}(B)=\mathsf{H}_{12}(B,B)=\sum^3_{i=1}\mathsf{H}_{12}(B_i)+\sum_{i\neq j} \mathsf{H}_{12}(B_i,B_j)=0,
\]
where ``cross-helicities'' $\mathsf{H}_{12}(B_i,B_j)$, $i\neq j$ vanish by Formula \eqref{eq:helicity12-flux-tubes}, and self helicities $\mathsf{H}_{12}(B_i)$ vanish because the average linking number of orbits is zero (orbits are just parallel 
circles). Nevertheless, Formula \eqref{eq:helicity-fluxes} tells us 
\[
\mathsf{H}_{123}(B;\mathcal{T}^\text{Borr})\neq 0,
\]
as a result we may regard $\mathsf{H}_{123}(B;\mathcal{T})$ as a possible ``higher obstruction'' to the energy relaxation or the third order cross-helicity of $B$ on $\mathcal{T}$.

 To obtain a lower bound for $E_2(B)$ in such situations we notice that 
$\mathsf{H}_{123}(B;\mathcal{T})$ is the $L^2$-inner product of the $6$-forms:
$\ast(\omega\wedge\alpha_\omega)$ and $\iota_{B_1}\mu_1\wedge \iota_{B_2}\mu_2\wedge \iota_{B_3}\mu_3$.
Indeed, from \eqref{eq:mu_123-helicity} we have
\[
 \mathsf{H}_{123}(B;\mathcal{T})=\bigl(\ast(\alpha_\omega\wedge \omega),\iota_{B_1}\mu_1\wedge \iota_{B_2}\mu_2\wedge \iota_{B_3}\mu_3\bigr)_{L^2\Omega^6((S^3)^3)}\ .
\]
Let $C_\mathcal{T}=\|\alpha_\omega\wedge \omega\|_{L^2}$, for a fixed Riemannian product metric on $(S^3)^3$ this constant depends on the domain $\mathcal{T}$ in $(S^3)^3$. We estimate using the Cauchy-Schwarz inequality
\begin{equation}\label{eq:energy-est}
\begin{split}
 |\mathsf{H}_{123}(B;\mathcal{T})| & \leq  \|\alpha_\omega\wedge \omega\|_{L^2}\|\iota_{B_1}\mu_1\wedge \iota_{B_2}\mu_2\wedge \iota_{B_3}\mu_3\|_{L^2}\\
 & =  C_\mathcal{T}\Bigl(\int_{\mathcal{T}}(\iota_{B_1}\mu_1\wedge \iota_{B_2}\mu_2\wedge \iota_{B_3}\mu_3)\wedge \ast(\iota_{B_1}\mu_1\wedge \iota_{B_2}\mu_2\wedge \iota_{B_3}\mu_3)\Bigr)^\frac{1}{2}\\
 & \stackrel{(1)}{=}  C_\mathcal{T}\Bigl(\int_{\mathcal{T}}|B_1|^2|B_2|^2|B_3|^2\,\mu\Bigr)^\frac{1}{2}\\
 & \stackrel{(2)}{=} C_\mathcal{T}\Bigl(\int_{\mathcal{T}_1}|B_1|^2\Bigr)^\frac{1}{2}\Bigl(\int_{\mathcal{T}_2}|B_2|^2\Bigr)^\frac{1}{2}\Bigl(\int_{\mathcal{T}_3}|B_3|^2\Bigr)^\frac{1}{2}\\
 & \leq    C_\mathcal{T} E_2(B)^\frac{3}{2},
\end{split}
\end{equation}
where $E_2(B)=\int_{S^3} |B|^2$.
To observe (1) first note that for any pair of forms: $\alpha\in \Omega^k(M)$, and $\beta\in \Omega^j(N)$, on Riemannian manifolds $M$ and $N$, on the product $M\times N$ we have  \begin{equation}\label{eq:hodge-product}
  \ast_{M\times N}(\pi^\ast_M(\alpha)\wedge \pi^\ast_N(\beta))=(\ast_M \pi^\ast_M \alpha)\wedge (\ast_N \pi^\ast_N \beta),
 \end{equation}
where $\pi_M:M\times N\longrightarrow M$ and $\pi_N:M\times N\longrightarrow N$ are the natural projections, (the proof is a simple calculation in an orthogonal frame of the product and is left to the reader).
Now, step (1) in \eqref{eq:energy-est} follows by applying \eqref{eq:hodge-product} to the integrand, and observing in the coframe $\{\eta^i_k\}$:
\begin{eqnarray*}
 \iota_{B_i}\mu_i\wedge \ast \iota_{B_i}\mu_i & = & (a_1\eta^i_2\wedge\eta^i_3-a_2\eta^i_1\wedge\eta^i_3+a_3\eta^i_1\wedge\eta^i_2)\wedge (a_1\eta^i_1-a_2\eta^i_2+a_3\eta^i_3)\\
 & = & (a^2_1+a^2_2+a^2_3)\, \eta^i_1\wedge\eta^i_2\wedge\eta^i_3 = |B_i |^2 \eta^i_1\wedge\eta^i_2\wedge\eta^i_3=|B_i |^2 \mu_i,
\end{eqnarray*}
where  $B_i=(a_1,a_2,a_3)$. Step (2) in \eqref{eq:energy-est} follows from Fubini Theorem.

Next, we aim to provide an estimate for $C_\mathcal{T}$. For this purpose we review 
some basic $L^2$-theory of the operator $d^{-1}$ (i.e. inverse of the exterior derivative $d$).
The main goal is to estimate an $L^2$-norm of the potential $\alpha_\omega$ of $\omega$ in \eqref{eq:omega-exact}. Following the standard elliptic theory of differential forms, \cite{Schwarz95}, the potential $\alpha_\omega$ in \eqref{eq:omega-exact} can be obtained via a solution to the Neumann problem for $2$-forms on $\widetilde{\mathcal{T}}$ (see Appendix \ref{apx:C})
\begin{equation}\label{eq:neumann-problem}
 \begin{cases}
  \Delta \phi_N=\omega,\qquad \text{in } \widetilde{\mathcal{T}}\\
  \mathbf{n}\,\phi_N=\mathbf{n}\,d\phi_N=0,\qquad \text{on } \partial \widetilde{\mathcal{T}},
 \end{cases}
\end{equation}
where $\mathbf{n}$ stands for the normal component of a differential form, and $\Delta=d\delta+\delta d$, $\delta=\pm\ast d\ast$, (c.f. \cite{Schwarz95}).
As $\mathcal{T}$ is a domain with corners we replace it by a slightly larger domain $\widetilde{\mathcal{T}}$ in $(S^3)^3$ with the same topology 
(i.e. $\mathcal{T}$ is a deformation retract of $\widetilde{\mathcal{T}}$) but with smooth boundary $\partial\widetilde{\mathcal{T}}$. (One may argue that it is not really necessary, since $\mathcal{T}$ is 
Lipschitz and elliptic problems, such as \eqref{eq:neumann-problem} are well posed on Lipschitz domains, \cite{Mitrea-Mitrea-Taylor01}).
Because of  $(i)$ in Theorem \ref{th:invariance-th}, we may use the restriction of 
\begin{equation}\label{eq:potential-neumann}
\alpha_\omega=\delta\phi_N,
\end{equation}
 to $\mathcal{T}$ (see Appendix \ref{apx:C} for justification of \eqref{eq:potential-neumann}).
Associated to \eqref{eq:neumann-problem} is the Neumann Laplacian
\[
 \Delta_N:H^2\Omega^2_N(\widetilde{\mathcal{T}})\longrightarrow L^2\Omega^2(\widetilde{\mathcal{T}}),
\] 
which has a discrete positive spectrum $\{\lambda_{i,N}\}$ and eigenvalues satisfy the variational principle called \emph{Rayleigh-Ritz quotient}, \cite{Chavel84}. The first (principal) eigenvalue  $\lambda_{1,N}$ may be expressed as
\[
 \lambda_{1,N}=\inf\left( 	\frac{\|d\varphi\|^2_{L^2}+\|\delta\varphi\|^2_{L^2}}{\|\varphi\|^2_{L^2}}\ \Bigl|\ \varphi\in H^1\Omega^2_N(M)\cap \mathcal{H}^2_N(M)^\perp\right)
\]
\no We denote the inverse of $\Delta_N$ by
\[
 G_N:L^2\Omega^2(\widetilde{\mathcal{T}})\longrightarrow H^2\Omega^2_N(\widetilde{\mathcal{T}}),
\]
which restricts to a compact, self-adjoint operator on $L^2$. As a result the spectrum of $G_N$ is 
discrete and given as $\{1/\lambda_{i,N}\}$. Note that based on these considerations we may define $d^{-1}:=\delta G_N$.

\begin{theorem}[Energy bound]
 For every volume preserving vector field $B$ which has an invariant unlinked domain $\mathcal{T}$,
 the $L^2$-energy of $B$ on $S^3$ is bounded below by the third order helicity $\mathsf{H}_{123}(B;\mathcal{T})$, as follows
 \begin{equation}\label{eq:H123-estimate}
  |\mathsf{H}_{123}(B;\mathcal{T})|\leq C_\mathcal{T} (E_2(B))^\frac{3}{2}.
 \end{equation}
Also, we may estimate the constant $C_\mathcal{T}$:
 \begin{equation}\label{eq:C_T-estimate}
  C_\mathcal{T}\leq \frac{1}{\sqrt{\lambda_{1,N}}}\|\omega\|^2_{L^\infty\Omega^2(\widetilde{\mathcal{T}})}(\text{Vol}(S^3))^3,
 \end{equation}
where $\lambda_{1,N}$ is the first eigenvalue of the Neumann  Laplacian on  $\Omega^2(\widetilde{\mathcal{T}})$.
\end{theorem}
\begin{proof}
We estimate
\begin{eqnarray*}
 \|\alpha_\omega\|^2_{L^2} = \|\delta \phi_N\|^2_{L^2} & \leq & \|d\phi_N\|^2_{L^2}+\|\delta\phi_N\|^2_{L^2}\\
 & = & |(\Delta \phi_N,\phi_N)| \leq  \|\omega\|_{L^2}\|\phi_N\|_{L^2},
\end{eqnarray*}
where we used the Green's formula \cite[p. 60]{Schwarz95} and boundary conditions of \eqref{eq:neumann-problem} in the second identity. Now, because $\phi_N=G_N\omega$, and 
it is a well known fact that $\|G_N\|_{L^2}=\frac{1}{\lambda_{1,N}}$, ($G_N$ is compact self adjoint on $L^2$) we obtain:
\[
\|\alpha_\omega\|_{L^2}\leq \frac{1}{\sqrt{\lambda_{1, N}}} \|\omega\|_{L^2}\ .
\]
As a result we estimate $C_\mathcal{T}$:
\begin{eqnarray*}
 C_\mathcal{T}=\|\alpha_\omega\wedge\omega\|_{L^2} & \leq & \|\omega\|_{L^\infty}\|\alpha_\omega\|_{L^2}\\
 & \leq &  \frac{1}{\sqrt{\lambda_{1, N}}} \|\omega\|_{L^2}\|\omega\|_{L^\infty}.
\end{eqnarray*}
\end{proof}
Notably, the best energy estimate so far has been obtained by Freedman and He, \cite{Freedman91-2}, for the $L^{3/2}$-energy of $B$,  in the case when $B$ admits an invariant domain $\mathcal{T}$ modeled on an $n$-component link $L=\{L_1,\ldots,L_n\}$ in $\R^3$. Their estimate
is based on the asymptotic crossing number and reads
\begin{equation}\label{eq:energy-free-he}
 E_{3/2}(B)\geq \Bigl(\frac{\pi}{16}\Bigr)^{1/4}\left(\sum^n_{k=1} \text{ac}(L_k,L)|\text{Flux}(B_k)|\right)\cdot \min_{1\leq k\leq n}\{|\text{Flux}(B_k)|\},
\end{equation}
where the asymptotic crossing numbers $\text{ac}(L_k,L)$  for Borromean links can be estimated below by
a smallest genus among surfaces in $\R^3\setminus \{L_1\cup\ldots\cup\hat{L}_k\cup\ldots\cup L_n\}$
with a single boundary component $L_i$.  Since $L^{3/2}$-energy of $B$ bounds the $L^2$-energy, inequality
\eqref{eq:energy-free-he} leads to a lower estimate purely in terms of fluxes and topological data.
It is not clear to the author if this approach can be extended to the case of invariant handlebodies 
considered in Section \ref{sec:handlebody}.
A different, more optimal estimate, has been obtained by Laurence and Stredulinsky, via the Massey product formula, in \cite{Laurence-Stredulinsky00a}, but the proof is provided only in a special case of the vector field $B$.

Contrary to these lower bounds, which are given in terms of topological data, the estimate in \eqref{eq:H123-estimate} depends on the geometry of the domain $\mathcal{T}$, and also $\|\omega\|_{L^\infty}$. Unfortunately, $\omega$ blows up on the diagonals $\mathbf{\Delta}\subset (S^3)^3$, and as a result the estimate is meaningless when the handlebodies 
$\mathcal{T}_i$ get close to each other during the evolution of the magnetic field $B$. At this point, we need an assumption for $\mathcal{T}_i$ to stay \emph{1cm} apart during the evolution.
Another drawback is that $\lambda_{1,N}$ is a geometric constant which is altered during the evolution as well. If we consider the situation in which the boundaries of $\mathcal{T}_i$ are invariant during the evolution, the estimate may be useful. Under such assumption, which occurs whenever the velocity field $v$ of plasma in \eqref{eq:Eulers-eq} is tangent to $\partial\mathcal{T}_i$, the bound in \eqref{eq:C_T-estimate} stays constant.

\section{Comparison to the known approaches via Massey products.}\label{sec:massey}

In several prior works \cite{Berger90, Berger-Evans92, Laurence-Stredulinsky00b, Mayer03} helicities were developed via the Massey product formula for $\bar{\mu}_{123}$. These approaches are equivalent to the one presented here in the sense that invariants obtained this way measure the same topological information. Most notably the work \cite{Laurence-Stredulinsky00b} provides an explicit expression for the third order helicity of the Borromean flux tubes, where the ergodic interpretation in the style of Arnold's asymptotic linking number is also provided.  In  \cite{Mayer03} one finds the following formula for the third order helicity
\begin{equation}\label{eq:massey-3A}
 \mathsf{M}_{123}(B_1,B_2,B_3)=\int_M A_1\wedge A_2\wedge A_3,
\end{equation}
where $A_i=d^{-1}(\iota_{B_i}\mu)$. This formula is valid for three distinct vector fields $B_i$
on a closed manifold $M$. For invariant domains with boundary 
\eqref{eq:massey-3A} defines an invariant provided $A_i\bigl|_{\partial M}=0$, but this only happens in certain situations (e.g. $M$ is simply connected, and $A_i$'s are appropriately chosen). 

\no The most commonly known formula directly related to the Massey products was developed by Berger \cite{Berger90}, in the case of 
Borromean flux tubes $\mathcal{T}=\mathcal{T}_1\cup \mathcal{T}_2\cup \mathcal{T}_3$
\begin{equation}\label{eq:massey123}
 \mathsf{M}_{123}(B_1,B_2,B_3)=\int_{\partial \mathcal{T}_1} A_1\wedge F_{23}+F_{12}\wedge A_3.
\end{equation}
In \cite{Berger90} it is expressed as a volume integral over $\mathcal{T}$ by applying \emph{gauge fixing}. When $\mathcal{T}_i$ are topologically solid tori there exists a single Massey product $<a_1,a_2,a_3>$ in the complement $S^3\setminus \mathcal{T}$, represented by the 2-form $A_1\wedge F_{23}+F_{12}\wedge A_3$. When $\mathcal{T}_i$ are handlebodies there are multiple Massey products, but the formula 
\eqref{eq:massey123} should still be valid. So far, such extensions have not been considered in the literature and the volume integrals over $\mathcal{T}$ may be harder to obtain in such a case.
One may also point out that ergodic interpretations 
of Massey products are more involved \cite{Laurence-Stredulinsky00b} comparing to the approach presented in Section \ref{sec:ergodic}.



\section*{Appendices}

\appendix


\section{Equations for the \emph{frozen-in-field} forms}\label{apx:A}
Given a volume preserving vector field $B$ on $M$ and a path $t\longrightarrow g(t)\in \text{\rm SDiff}(M)$, let 
\[
 B^t:=g_\ast(t)B\ . 
\]
Then for $B^0=B$, by definition 
\begin{equation}\label{eq:euler2}
 \dot{B^t}=\frac{d}{dt} B^t\bigl|_{t=0}=-\Lie_V B=-[V,B],
\end{equation}
and as a result 
\[
  \frac{d}{dt} \iota_{B^t}\mu|_{t=0}  =  -\iota_{[V,B]} \mu  
\]
Next, calculate ($\partial_t+\Lie_V\equiv \Lie_{\partial_t+V}$)
\begin{eqnarray*}
 (\partial_t+\Lie_V)\iota_{B^t}\mu & = & \partial_t(\iota_{B^t}\mu)+\Lie_V\iota_{B^t}\mu\\
 & = & \iota_{\dot{B^t}}\mu+\iota_{B^t}(\Lie_V\mu)+\iota_{[V,B^t]}\mu\\
 & = & \iota_{\dot{B^t}+[V,B^t]}\mu,
\end{eqnarray*}
where in the second identity we applied the general formula: $\iota_{[A,B]}=\Lie_A\iota_B-\iota_B\Lie_A$, and in the third equation the fact that $V$ is volume preserving i.e. $\Lie_V\mu=0$. As a result 
of \eqref{eq:euler2} we obtain
\[
 (\partial_t+\Lie_V)\iota_{B^t}\mu=0\ .
\]

\no Next, we justify Formula \eqref{eq:iota-volume}. First use \eqref{eq:iota-antiderivation} to calculate
\begin{equation}\label{eq:iota-vol1}
 \begin{split}
  \iota_{B_i} (\mu_j\wedge\mu_k) & = (\iota_{B_i} \mu_j)\wedge\mu_k-\mu_j\wedge(\iota_{B_i}\mu_k)\\
  \iota_{B_i} (\mu_1\wedge\mu_2\wedge\mu_3) & = (\iota_{B_i} \mu_1)\wedge\mu_2\wedge\mu_3-\mu_1\wedge \iota_{B_i}(\mu_2\wedge\mu_3)\\
  & =(\iota_{B_i} \mu_1)\wedge\mu_2\wedge\mu_3-\mu_1\wedge \iota_{B_i}\mu_2\wedge\mu_3+\mu_1\wedge \mu_2\wedge\iota_{B_i}\mu_3,
 \end{split}
\end{equation}
since $\iota_{B_i}\mu_j=0$ for $i\neq j$ only one term in the above expressions remains for each $i$.
Set $\alpha:=\iota_{B_2}\iota_{B_1}\beta$, $\beta\in \Omega^3((S^3)^3)$, since $\alpha\wedge \mu_1\wedge\mu_2\wedge\mu_3=0$ and $\alpha$ is a 1-form we obtain
\[
\begin{split}
 0 & =\iota_{B_3}(\alpha\wedge \mu_1\wedge\mu_2\wedge\mu_3)\\
   & =(\iota_{B_3}\alpha)\wedge\mu_1\wedge\mu_2\wedge\mu_3-\alpha\wedge \mu_1\wedge\mu_2\wedge\iota_{B_3}\mu_3,
\end{split}
\]
where in the last step we used \eqref{eq:iota-vol1}. Therefore
\[
(\iota_{B_3}\iota_{B_2}\iota_{B_1}\beta)\wedge\mu_1\wedge\mu_2\wedge\mu_3=(\iota_{B_2}\iota_{B_1}\beta)\wedge \mu_1\wedge\mu_2\wedge\iota_{B_3}\mu_3\ .
\]
Analogously, $(\iota_{B_2}\iota_{B_1}\beta)\wedge \mu_1\wedge\mu_2=(\iota_{B_1}\beta)\wedge \mu_1\wedge \iota_{B_2}\mu_2$ and $\iota_{B_1}\beta\wedge \mu_1=\beta\wedge \iota_{B_1}\mu_1$ which justifies Equation \eqref{eq:iota-volume}.

\section{Zero contribution of short paths to the time average}\label{apx:B}

It is clear when $\beta\in \Omega^3(\widetilde{\mathcal{T}})$ is at least a $C^1$ on
$\widetilde{\mathcal{T}}\subset S^3\times S^3\times S^3$, then
$f=\beta(B_1,B_2,B_3)$ is continuous on $\widetilde{\mathcal{T}}$ and
\begin{gather*}
\int_{\bar{\mathscr{O}}^1_T(x)}\int_{\bar{\mathscr{O}}^2_T(y)}\int_{\bar{\mathscr{O}}^3_T(z)}
f=\\
\bigl(\int_{\mathscr{O}^1_T(x)}+\int_{\sigma(\Phi^1(x,T),x)}\bigr)
\bigl(\int_{\mathscr{O}^2_T(y)}+\int_{\sigma(\Phi^2(y,T),y)}\bigr)
\bigl(\int_{\mathscr{O}^3_T(z)}+\int_{\sigma(\Phi^3(z,T),z)}\bigr)f\ .
\end{gather*}
After expanding, it is obvious that we must show the following (for all choices), when $T\to \infty$:
\begin{eqnarray}\label{eq:lim-sp1}
\frac{1}{T^3}\int_{\sigma(\Phi^1(x,T),x)}\int_{\mathscr{O}^2_{T}(y)}\int_{\mathscr{O}^3_{T}(z)} f & \longrightarrow & 0,\\
\label{eq:lim-sp2}\frac{1}{T^3}\int_{\sigma(\Phi^1(x,T),x)}\int_{\sigma(\Phi^2(y,T),y)}\int_{\mathscr{O}^3_{T}(z)} f & \longrightarrow & 0,\\
\label{eq:lim-sp3}\frac{1}{T^3}\int_{\sigma(\Phi^1(x,T),x)}\int_{\sigma(\Phi^2(y,T),y)}\int_{\sigma(\Phi^3(z,T),z)} f & \longrightarrow & 0.
\end{eqnarray}
Since the lengths of the short paths in $\mathcal{S}$ are bounded by a common constant $d$, \eqref{eq:lim-sp1}-\eqref{eq:lim-sp3} follow immediately, e.g. for \eqref{eq:lim-sp1} we have

\begin{gather*}
\bigl|\frac{1}{T^3}\int_{\sigma(\Phi^1(x,T),x)}\int_{\mathscr{O}^2_{T}(y)}\int_{\mathscr{O}^3_{T}(z)} f\bigr|\leq \frac{1}{T^3} d(T+d)(T+d)\|f\|_{\infty} \longrightarrow  0.
\end{gather*}

\section{Notation in Section \ref{sec:energy}}\label{apx:C}
 We adopt notation from the elegant exposition in \cite{Schwarz95}. Let $M$ be an orientable manifold with smooth boundary.
\begin{itemize}
\item[] $\Omega^k(M)=C^\infty(M,\Lambda^k)$, smooth differential forms on $M$.
\item[] $\Omega^k_N(M)=\{\phi\in \Omega^k(M)\ |\ \mathbf{n}\phi=0,\, \mathbf{n}d\phi=0\}$ the subspace satisfying the Neumann boundary conditions, ($\mathbf{n}$ denotes a normal component of a form along $\partial M$).
\end{itemize}
\no The $L^2$-inner product on $\Omega^k(M)$ is defined as 
 \[
  \bigl( \omega,\eta\bigr)_{L^2}=\int_M \omega\wedge\ast \eta,
 \] 
\begin{itemize}
\item[] $L^2\Omega^k(M)$, $L^2$-differential forms on $M$.
\item[] $\mathcal{H}^k_N(M)=\{\lambda\in H^1\Omega^k(M)\ |\ d\lambda=\delta\lambda=0, \mathbf{n}\lambda=0\}$ the subspace of the Neumann harmonic fields.
\end{itemize}

\no Next we justify \eqref{eq:potential-neumann}, first observe that for any $\gamma\in H^1\Omega^{k-1}(M)$:
\begin{equation}\label{eq:perp-neumann}
 (d\gamma,\lambda)_{L^2}=\int_{\partial M} \mathbf{t}\gamma\wedge\ast\mathbf{n}\lambda=0,\qquad \forall \lambda\in \mathcal{H}^k_N(M),
\end{equation}
where $\mathbf{t}$, and $\mathbf{n}$ stands for respectively tangent and normal to $\partial M$ components of the form. As a result, if $\omega \in \mathcal{H}^k_N(M)^\perp$ we obtain a solution $\phi$ to the Neumann problem:
\begin{equation}\label{eq:neumann-phi}
 \delta d\phi+d\delta\phi=\omega,\quad \Rightarrow\quad \omega-d\delta\phi=\delta d\phi\ .
\end{equation}
Formula \eqref{eq:perp-neumann} implies: $(\omega-d\delta\phi)\in \mathcal{H}^k_N(M)^\perp$, moreover
$\mathbf{n}(\omega-d\delta\phi)=\mathbf{n}(\delta d\phi)=\delta\,\mathbf{n}(d\phi)=0$, by the boundary condition in \eqref{eq:neumann-problem}. If $\omega$ is a closed form, $\omega-d\delta\phi$ is also closed,
and clearly coclosed by \eqref{eq:neumann-phi}. Thus $\omega-d\delta\phi$ is a harmonic field with zero normal component, and therefore it has to be in $\mathcal{H}^k_N(M)$, and therefore the zero form. This yields 
\[
 \omega=d\delta\phi\ .
\]
As a result we obtain a necessary and sufficient conditions for $\omega$ to be exact:
 \begin{itemize}
 \item[$(i)$] $d\omega=0$,
 \item[$(ii)$] $(\omega,\lambda)_{L^2}=0$, for all $\lambda\in \mathcal{H}^k_N(M)$.
\end{itemize}

\vspace{.5cm}

\no\underline{\hspace{.4\textwidth}}

\no {\sc Department of Mathematics,\ Univ. of Pennsylvania,\ Philadelphia,\  PA, 19104}\\
e-mail: \texttt{rako@math.upenn.edu}\\
URL: \texttt{www.math.upenn.edu/\~{}rako}

\end{document}